\documentclass[preprint,3p,12pt]{elsarticle}
\usepackage{mathrsfs}
\usepackage{amsmath}
\usepackage{stmaryrd}
\usepackage{bbding}
\usepackage{dcolumn}
\usepackage{graphicx}
\usepackage{amsfonts}
\usepackage{amssymb}
\usepackage{psfrag}
\usepackage{wrapfig}
\usepackage{subfigure}
\usepackage{makeidx}
\usepackage{bm}
\usepackage{epsf}
\usepackage{epsfig}
\usepackage{setspace}
\usepackage{graphicx}
\usepackage{epstopdf}
\usepackage{psfrag}
\usepackage{subfigure}
\usepackage{color}

\begin{document}
\title{A Well-Balanced Space-Time ALE Compact Gas-Kinetic Scheme for the Shallow Water Equations on Unstructured Meshes}

\author[HKUST1]{Fengxiang Zhao}
\ead{fzhaoac@connect.ust.hk}

\author[HKUST1,HKUST2]{Jianping Gan}
\ead{magan@ust.hk}

\author[HKUST1,HKUST2]{Kun Xu\corref{cor}}
\ead{makxu@ust.hk}

\address[HKUST1]{Department of Mathematics, Hong Kong University of Science and Technology, Clear Water Bay, Kowloon, Hong Kong}
\address[HKUST2]{Center for Ocean Research in Hong Kong and Macau (CORE), Hong Kong University of Science and Technology, Clear Water Bay, Kowloon, Hong Kong}
\cortext[cor]{Corresponding author}

\begin{abstract}
This study presents a high-order, space-time coupled arbitrary Lagrangian Eulerian (ALE) compact gas-kinetic scheme (GKS) for the shallow water equations on moving unstructured meshes. The proposed method preserves both the geometric conservation law (GCL) and the well-balanced property.
Mesh motion effects are directly incorporated by formulating numerical fluxes that account for the spatial temporal nonuniformity of the flow field and the swept area of moving cell interfaces.
This allows temporal updates to be performed on the physical moving mesh, avoiding data remapping.
The compact GKS provides time accurate evolution of flow variables and fluxes, enabling the scheme to achieve second-order temporal accuracy within a single stage.
To consistently treat bottom topography on moving meshes, an evolution equation for the topography is established and discretized using a compatible space-time scheme, in which the fluxes induced by mesh motion are computed accurately. Mathematical proofs demonstrating the GCL preserving and well-balanced properties of the proposed ALE formulation are also provided. For improved accuracy and robustness, a nonlinear fourth-order compact reconstruction technique is employed. A comprehensive set of numerical experiments verifies the scheme's theoretical properties and demonstrates its accuracy, stability, and effectiveness in simulating complex shallow-water flow problems.

\end{abstract}


\maketitle

\section{Introduction}

The shallow water equations (SWE) provide a canonical model for geophysical water flows. Applications range from large-scale ocean circulation to regional coastal and channel hydraulics, including tsunamis, tides, storm surges, and dam-break waves. Accordingly, a broad suite of numerical schemes has been developed for solving the SWE \cite{leveque-1998,schwanenberg2000dg,zhou2001-surface,xu2002-swe,xing2005_highFD}.
For realistic engineering and geophysical simulations, the standard SWE must often be augmented with models for additional physical processes. These may include flow stratification \cite{castro2001_TL,abgrall2009_TL}, sediment transport and morphodynamics \cite{liu2017coupled,garcia2019_friction}, fluid-structure interaction \cite{de2014_FSI,haidar2024_moving}, and external forcings such as wind stress, the Coriolis effect, and bottom friction \cite{mandli2014_storm-surge}. Consequently, the development of numerical methods capable of robustly handling these coupled complexities is crucial for advancing the predictive capabilities of SWE models.

A central algorithmic requirement is the well-balanced property. The property requires a precise discrete balance between numerical fluxes and source-term contributions that preserves steady states (e.g., the lake-at-rest solution). Additional challenges include sharply capturing hydraulic jumps and other discontinuities, resolving fine-scale features, and robustly handling moving wet-dry fronts, among others. High accuracy is further demanded by the multiscale nature of geophysical flows, where characteristic wavelengths of interest may be far smaller than computational domains spanning hundreds of kilometers. Effective strategies for efficient, accurate simulation therefore include high-order schemes \cite{brus2019high_application}, adaptive moving-mesh methods \cite{porta2012anisotropic}, and adaptive mesh refinement techniques \cite{mandli2014_storm-surge}.

In this study, we develop a well-balanced, high-order space-time ALE scheme for the SWE on moving unstructured meshes, suitable for moving-boundary problems and adaptive high-resolution simulations. On moving meshes, numerical schemes for hyperbolic conservation laws must satisfy the GCL, which requires that, under a uniform free stream, the flux induced solely by mesh motion exactly equals the face-swept measure. For a flat bottom, the source term vanishes and the SWE reduce to a system of hyperbolic conservation laws; thus, any SWE scheme must satisfy the GCL. With non-flat bathymetry, however, the topography evolves relative to the moving mesh and mesh motion complicates both flux and source contributions, making the construction of a well-balanced discretization substantially more challenging.
Within the compact GKS framework, we design a high-order space-time ALE method that preserves both the GCL and the well-balanced property. The time-dependent gas distribution function at cell interfaces provides time-accurate numerical fluxes and interface states, while its consistently derived spatial gradients encode spatial nonuniformity.
These quantities are uniformly obtained from the distribution function, without additional algorithmic complexity or cost.
Leveraging these quantities enables an accurate treatment of mesh-motion effects and associated source terms, yielding a discretization that preserves both the GCL and the well-balanced property. Owing to the space-time ALE formulation, the proposed compact GKS updates the solution directly on the physical moving mesh without data remapping. The resulting space-time moving-mesh methodology is broadly applicable to general hyperbolic conservation laws.

The second-order GKS was first developed for compressible flow simulations \cite{xu2,xu1}.
The GKS provides a unified framework for computing both inviscid and viscous fluxes, features a time-dependent and intrinsic multi-scale mechanism to handle smooth flows and discontinuities alike, and achieves second-order accuracy in a single time-stepping procedure based on the gas distribution function.
Furthermore, the GKS's core principle, the time-accurate evolution of the gas distribution function at cell interfaces, has enabled the recent development of high-order compact schemes, demonstrating its versatility and potential for high-fidelity simulations \cite{CGKSAIA,zhao2023direct}.
For example, the compact GKS from fourth- to sixth-order GKS on triangular mesh was proposed in \cite{zhaocompact_tri}.
In the compact GKS, both cell-averaged flow variables and their gradients are updated by the conservation laws and Gauss's theorem, which helps the compact spatial reconstruction with involving only neighboring cells.

The compact GKS for compressible flows has been extended to the SWE in \cite{zhao2021-swe}, where complex models representing realistic effects are incorporated as source terms. The scheme employs the same fourth-order compact reconstruction and fourth-order temporal discretization as in \cite{zhaocompact_tri}.
Because the GKS flux is time-dependent and already incorporates source effects, designing well-balanced GKS schemes is particularly challenging.
In \cite{zhao2024_TLSWE}, the fourth-order compact scheme is also extended to solve the two-layer SWE. Compared to the single-layer SWE, the two-layer SWE introduces more complex source terms, posing challenges in constructing the numerical scheme.
As a representative method among the multistage and multiderivative methods, the 2-stage 4th-order (S2O4) method \cite{li2019AIA} is adopted in the compact GKS for achieving high-order temporal accuracy with less stages \cite{zhaocompact_tri,zhao2021-swe}.

This paper is organized as follows. Section 2 presents the governing equations and the space-time ALE discretization of the SWE on moving meshes.
Section 3 is about the GCL and well-balanced property of the scheme.
The high-order compact spatial reconstruction will be presented in Section 4.
Section 5 is about the adaptive mesh moving method.
In Section 6, a series of numerical examples will validate the ALE compact GKS on moving unstructured meshes.
Section 7 is the conclusion.

\section{Space-time coupled numerical formulations on a moving mesh}

This section presents a GCL-preserving, well-balanced compact GKS for the SWE on moving meshes, formulated within a space-time coupled ALE framework.
The method is built on the gas distribution function of the GKS, which provides time-accurate interface fluxes and macroscopic states at cell interfaces.
Unlike conventional schemes that decouple space and time, the proposed space-time coupled compact GKS updates the solution directly on the physical moving mesh without data remapping, while accurately capturing flow unsteadiness and spatial nonuniformity.
For non-flat bathymetry, the bed elevation is advanced consistently with the mesh motion by introducing a linear advection equation for the topography and discretizing it within the same space-time framework.

\subsection{SWE on a moving mesh}

The SWE, spatially integrated over each mesh cell, are given by
\begin{equation}\label{SWE}
\int_{\Omega_{(t)}} \frac{\partial \textbf{W}}{\partial t} \text{d}\Omega= -\int_{\partial\Omega_{(t)}} \textbf{F} \cdot \textbf{n} \text{d}s +\int_{\Omega_{(t)}} \textbf{S} \text{d}\Omega,
\end{equation}
where $\partial\Omega_{(t)}$ is the interface of the moving mesh cell $\Omega_{(t)}$.
$\mathbf{W}$ denotes the vector of water height and momentum, and $\mathbf{F}=(\textbf{F}^x,\textbf{F}^y)$ is the corresponding flux, with $\textbf{F}^x$ and $\textbf{F}^y$ representing the fluxes in the x and y directions, respectively.
The source term $\mathbf{S}$ accounts for the forcing induced by non-flat bottom topography. In cases where additional physical effects are considered, $\mathbf{S}$ may also include other forcings \cite{zhao2021-swe}, such as bottom friction.
The specific forms of $\mathbf{W}$, $\mathbf{F}$, and $\mathbf{S}$ are
\begin{equation*}
{\textbf{W}} =
\left(
\begin{array}{c}
h\\
h U\\
h V\\
\end{array}
\right), \\
{\textbf{F}^x} =
\left(
\begin{array}{c}
hU\\
h U^2+\frac{1}{2}Gh^2\\
h UV\\
\end{array}
\right),\\
{\textbf{F}^y} =
\left(
\begin{array}{c}
hV\\
h UV\\
h V^2+\frac{1}{2}Gh^2\\
\end{array}
\right),
\end{equation*}
and
\begin{equation*}
{\textbf{S}} =
\left(
\begin{array}{c}
0\\
-GB_x h\\
-GB_y h\\
\end{array}
\right).
\end{equation*}

Taking into account that the cell $\Omega_{(t)}$ is moving, to facilitate the discretization of Eq. (\ref{SWE}), one approach is to reformulate it so that the time derivative on the left-hand side is moved outside the integral. This can be achieved via the Reynolds transport theorem, which states
\begin{equation}\label{Reynolds-transport-theorem}
\frac{\text{d}}{\text{d}t}\int_{\Omega_{(t)}}\textbf{W}\text{d}\Omega =\int_{\Omega_{(t)}} \frac{\partial \textbf{W}}{\partial t} \text{d}\Omega +\int_{\partial\Omega_{(t)}} \textbf{W} \textbf{V}^m\cdot \textbf{n} \text{d}s,
\end{equation}
where $\textbf{V}^m=(U^m,V^m)$ is the mesh moving velocity of $\partial\Omega_{(t)}$.
The second term on the right-hand side (RHS) of Eq. (\ref{Reynolds-transport-theorem}) accounts for the introduction and removal of space from $\Omega_{(t)}$ due to its moving interfaces, i.e., the flux originating from mesh motion, denoted as
\begin{align*}
\textbf{F}^m= \textbf{W} \textbf{V}^m.
\end{align*}
Consequently, combining Eq. (\ref{SWE}) and Eq. (\ref{Reynolds-transport-theorem}), the governing equations for the total derivative of the flow variables $\mathbf{W}$ over the moving cell $\Omega_{(t)}$ are derived as
\begin{equation}\label{SWE-moving-cell}
\frac{\text{d}}{\text{d}t}\int_{\Omega_{(t)}}\textbf{W}\text{d}\Omega =-\int_{\partial\Omega_{(t)}} \textbf{F} \cdot \textbf{n} \text{d}s +\int_{\partial\Omega_{(t)}} \textbf{F}^m \cdot \textbf{n} \text{d}s +\int_{\Omega_{(t)}} \textbf{S} \text{d}\Omega.
\end{equation}
For convenience, the RHS of Eq. (\ref{SWE-moving-cell}) is denoted as $\mathcal{L}$,
\begin{align*}
\mathcal{L}=-\int_{\partial\Omega_{(t)}} \textbf{F} \cdot \textbf{n} \text{d}s +\int_{\partial\Omega_{(t)}} \textbf{F}^m \cdot \textbf{n} \text{d}s +\int_{\Omega_{(t)}} \textbf{S} \text{d}\Omega.
\end{align*}

\begin{figure}[!htb]
\centering
\includegraphics[width=0.40\textwidth]{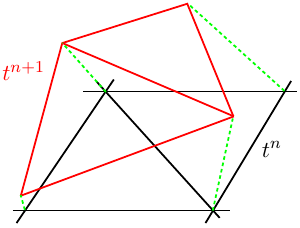}
\caption{\label{0-cell-moving} Schematic of the moving triangular mesh cell at times $t^n$ (black) and $t^{n+1}$ (red).}
\end{figure}

\subsection{Space-time ALE formulation}

This section presents the spatial and temporal discretizations of Eq. (\ref{SWE-moving-cell}).
The time discretization depends on the time-dependent numerical fluxes and flow variables obtained from the gas distribution function of GKS.
Integrating Eq. (\ref{SWE-moving-cell}) over time yields a one-stage, second-order time discretization, given by
\begin{equation}\label{SWE-moving-cell-time}
\textbf{W}_j^{n+1}|\Omega^{n+1}_j| =\textbf{W}_j^{n}|\Omega^{n}_j| +\mathcal{L}_j^n \Delta t + \frac{\mathrm{d} \mathcal{L}_j^n}{\mathrm{d} t}\frac{\Delta t^2}{2},
\end{equation}
where $\mathcal{L}_j^n$ is the discrete form of $\mathcal{L}$ on cell $\Omega_j$ at $t^n$.
Based on the one-stage second-order scheme and employing a S2O4 method \cite{zhao2021-swe}, a fourth-order temporal discretization is readily obtained.
The cell-averaged variable $\textbf{W}_j^{n}$ is defined as
\begin{align*}
\textbf{W}_j^{n}=\frac{1}{|\Omega^{n}_j|}\int_{\Omega^n_j}\textbf{W}(\mathbf{x},t^n)\text{d}\Omega.
\end{align*}
The line integrals in $\mathcal{L}$ can be approximated using high-accuracy Gaussian quadrature. In this study, a two-point Gaussian rule with fourth-order accuracy is employed.
However, since both the cell geometry and the flow variables evolve in time, evaluating the time derivative of $\mathcal{L}_j^n$ is nontrivial. This, in turn, makes an accurate implementation of the temporal discretization in Eq. (\ref{SWE-moving-cell-time}) particularly challenging.
The term $\mathcal{L}_j^n$ is discretized by considering its three constituent parts, denoted as $\mathcal{L}_j^n=\mathcal{L}_{1,j}^n+\mathcal{L}_{2,j}^n+\mathcal{L}_{3,j}^n$. The discrete formulation for each of these terms is detailed below.
For brevity, the superscript $n$ indicating time $t^n$ on the flow variables $\mathbf{W}$ and the fluxes $\mathbf{F}$ will no longer be shown.

\subsubsection{Hydrodynamic fluxes}
The first term $\mathcal{L}_{1,j}^n$ represents the contribution of the fluid convection itself to $\mathbf{W}_j$ over the control volume.
The discrete form of $\mathcal{L}_{1,j}^n$ is
\begin{equation}\label{SWE-discrete-L1}
\mathcal{L}_{1,j}^n = - \sum_{l=1}^{l_0}\big( \sum _{k=1}^2 \omega_k \textbf{F}(\textbf{x}_{l,k}) \cdot \textbf{n}_l\big) \big|\Gamma_{l} \big|,
\end{equation}
where $l_0$ is the number of faces of the cell (e.g., $l_0=3$ for a triangle), $|\Gamma_l|$ is the face length, $\mathbf{n}_l$ is the outward unit normal, and $\mathbf{x}_{l,k}$ and $\omega_k$ are the Gaussian quadrature points and weights (with $\omega_k=1/2$), respectively.
For notational convenience, hereafter, we will use $\sum_{l,k}(...)$ in the following text to represent $\sum_{l=1}^{l_0}\sum _{k=1}^2(...)$.
$\textbf{F}(\textbf{x}_{l,k})$ will be abbreviated as $\textbf{F}_{l,k}$. A similar convention is adopted for any other variable evaluated at $\mathbf{x}_{l,k}$.
The time derivative $\mathrm{d} \mathcal{L}_{1,j}^n/\mathrm{d} t$ can be obtained as
\begin{equation}\label{SWE-discrete-L1-dt}
\begin{split}
\frac{\mathrm{d} \mathcal{L}_{1,j}^n}{\mathrm{d} t}=& - \sum_{l,k}\big( \omega_k \partial \textbf{F}_{l,k}/\partial t \cdot \textbf{n}_l \big) \big|\Gamma_{l} \big| \\
&-\sum_{l,k}\big( \omega_k (\nabla P_{\mathrm{stat}})_{l,k} \mathbf{V}^m_{l,k} \cdot \mathbf{n}_l \big) \big|\Gamma_l\big| +Res(\mathcal{L}^n_{1,j}),
\end{split}
\end{equation}
where $\nabla P_{\mathrm{stat}}=(0,Ghh_x,Ghh_y)^\mathrm{T}$, and
\begin{align*}
\begin{split}
Res(\mathcal{L}^n_{1,j})=&-\sum_{l,k} \big( \omega_k \mathbf{V}^m_{l,k} \cdot \nabla\textbf{F}_{l,k} \cdot \textbf{n}_l \big) \big|\Gamma_{l} \big|
- \sum_{l,k} \omega_k \textbf{F}_{l,k} \cdot \mathrm{d} \big(\textbf{n}_l \big|\Gamma_{l} \big|\big)/\mathrm{d} t \\
&+\sum_{l,k}\big( \omega_k (\nabla P_{\mathrm{stat}})_{l,k} \mathbf{V}^m_{l,k} \cdot \mathbf{n}_l \big) \big|\Gamma_l\big|.
\end{split}
\end{align*}
Here $Res(\mathcal{L}^n_{1,j})$ denotes the contribution to $\mathrm{d} \mathcal{L}^n_{1,j}/\mathrm{d} t$ arising from changes in the direction and length of the moving interfaces, as well as from the spatial nonuniformity of the flow.
The numerical fluxes and their derivatives are based on the gas evolution solution of GKS (see Appendix).
When $\mathbf{F}(\mathbf{x})$ is constant over the cell, it can be shown, based on the closure condition of the cell, that $Res(\mathcal{L}^n_{1,j})=0$.
Leveraging this property and for algorithmic simplicity, we uniformly neglect $Res(\mathcal{L}^n_{1,j})$ in this study. As demonstrated in the numerical examples, omitting $Res(\mathcal{L}^n_{1,j})$ does not introduce noticeable numerical errors.

\begin{figure}[!htb]
\centering
\includegraphics[width=0.50\textwidth]{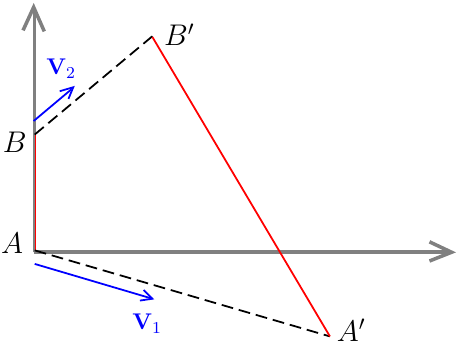}
\caption{\label{0-cell-moving} Schematic of the motion of edge $AB$ of a mesh cell. The motion is determined by the nodal velocities $V_1$ and $V_2$. Within a time step, the motion velocities are constant, and $A'B'$ remains a straight line.}
\end{figure}

\subsubsection{Mesh motion-induced fluxes}
The second term, $\mathcal{L}_{2,j}^n$, represents the flux generated by the change in cell area due to interface motion, where the motion involves translation and rotation. For the scope of this study, we consider interfaces that move with a constant velocity during each time step, while remaining straight, as depicted in Fig. \ref{0-cell-moving}.
The spatial discretization of $\mathcal{L}^n_{2,j}$ is given by
\begin{equation}\label{SWE-moving-cell-Fluxmove}
\mathcal{L}_{2,j}^n = \sum_{l,k}\big( \omega_k \textbf{F}^m_{l,k} \cdot \textbf{n}_l\big) \big|\Gamma_{l} \big|.
\end{equation}
The time derivative $\mathrm{d} \mathcal{L}_{2,j}^n/\mathrm{d} t$ is given by the chain rule as
\begin{equation}\label{SWE-moving-cell-Fluxmove-dt}
\begin{split}
\frac{\mathrm{d} \mathcal{L}_{2,j}^n}{\mathrm{d} t} = &\sum_{l,k}\big( \omega_k \frac{\mathrm{d}\textbf{W}_{l,k}}{\mathrm{d}t} \textbf{V}^m_{l,k} \cdot \textbf{n}_l \big) \big|\Gamma_{l} \big| \\
&+ \sum_{l,k} \omega_k \textbf{W}_{l,k} \frac{\mathrm{d}( \textbf{V}^m_{l,k} \cdot \textbf{n}_l |\Gamma_{l}| )}{\mathrm{d} t}.
\end{split}
\end{equation}
The first term on the RHS of Eq. (\ref{SWE-moving-cell-Fluxmove-dt}) arises from the unsteadiness of the flow and its spatial variation. The total derivative of $\mathbf{W}$ is given by
\begin{equation}\label{slope-W-evolution}
\frac{\mathrm{d} \mathbf{W}_{l,k}}{\mathrm{d}t}= \frac{\partial \mathbf{W}_{l,k}}{\partial t} + \textbf{V}^m_{l,k} \cdot \nabla\textbf{W}_{l,k},
\end{equation}
where $\frac{\partial \mathbf{W}_{l,k}}{\partial t}$ and $\nabla\textbf{W}_{l,k}$ at the interface are provided by the space-time coupled evolution solution of GKS, with details to be given in the Appendix.
Furthermore, the second term on the RHS of Eq. (\ref{SWE-moving-cell-Fluxmove-dt}) arises from the rotational and deformational effects of the interfaces motion.
For interface $\overrightarrow{A'B'}$ in Fig. \ref{0-cell-moving}, from the condition that the normal vector is always perpendicular to it, the unit outward normal vector at any time can be obtained as
\begin{align*}
\mathbf{n}_l=\frac{1}{|\Gamma_l|}(\Delta Vt +|\Gamma_l|,-\Delta Ut),
\end{align*}
where $\Delta U=U_2-U_1$ and $\Delta V=V_2-V_1$.
Then we can calculate
\begin{align*}
\begin{split}
\textbf{V}^m_{l,1} \cdot \textbf{n}_l |\Gamma_{l}|=& (\alpha \mathbf{V}_1+(1-\alpha)\mathbf{V}_2)\cdot \mathbf{n}_l |\Gamma_l| \\
=&(\alpha U_1+(1-\alpha)U_2)|\Gamma_l| +t (\mathbf{V}_1\times\mathbf{V}_2)\cdot \mathbf{k},\\
\textbf{V}^m_{l,2} \cdot \textbf{n}_l |\Gamma_{l}|=&(\alpha U_2+(1-\alpha) U_1)|\Gamma_l| +t(\mathbf{V}_1\times\mathbf{V}_2)\cdot \mathbf{k},
\end{split}
\end{align*}
where $\alpha=1/\sqrt 3$, $\mathbf{k}=\mathbf{i}\times \mathbf{j}$, and $\mathbf{i}$ and $\mathbf{j}$ are the basis vectors along two directions for the velocity vector.
The second term of the RHS in Eq. (\ref{SWE-moving-cell-Fluxmove-dt}) has been determined.
Thus the final discrete form of the time derivative of $\mathcal{L}_{2,j}^n$ is given by
\begin{equation}\label{SWE-moving-cell-Fluxmove-dt-1}
\begin{split}
\frac{\mathrm{d} \mathcal{L}_{2,j}^n}{\mathrm{d} t} =& \sum_{l,k}\big( \omega_k (\partial \textbf{W}_{l,k}/\partial t+\mathbf{V}^m_{l,k}\cdot \nabla\mathbf{W})_{l,k} \textbf{V}^m_{l,k} \cdot \textbf{n}_l \big) \big|\Gamma_{l} \big| \\
&+ \sum_{l,k} \omega_k \textbf{W}_{l,k} (\mathbf{V}_1\times\mathbf{V}_2)\cdot \mathbf{k}.
\end{split}
\end{equation}
It is worth noting that the interface values of $\mathbf{W}$ and their spatiotemporal derivatives can be uniformly obtained from the distribution function in GKS, without additional computational complexity or cost.

\subsubsection{Source terms}
Finally, the third term $\mathcal{L}_{3,j}^n$ in $\mathcal{L}_{j}^n$ is given by
\begin{equation}\label{SWE-moving-cell-S}
\mathcal{L}_{3,j}^n = -G h_j^n\big|\Omega_j^n \big| \big(0,B_{j,x}^n,B_{j,y}^n \big)^\mathrm{T},
\end{equation}
and
\begin{equation}\label{SWE-moving-cell-St}
\frac{\mathrm{d}\mathcal{L}_{3,j}^n}{\mathrm{d}t} = -G \frac{\mathrm{d}(h_j^n\big|\Omega_j^n \big|)}{\mathrm{d}t} \big(0,B_{j,x},B_{j,y} \big)^\mathrm{T}.
\end{equation}
The total derivative $\mathrm{d}\big(h_j^n\big|\Omega_j^n \big|\big)/\mathrm{d}t$ can be given exactly by the first component of Eq. (\ref{SWE-moving-cell-time}), as follows,
\begin{align*}
\frac{\mathrm{d}(h_j^n\big|\Omega_j^n \big|)}{\mathrm{d}t}=L_j^{n,h},
\end{align*}
where $L_j^{n,h}$ denotes the first component of $L_j^{n}$ corresponding to the RHS of updating $h_j$.

\subsection{Bottom topography update on a moving mesh}

For non-flat bottom topography, mesh motion causes the bottom distribution on the mesh unsteady, and the equation for the bottom elevation is given by
\begin{equation}\label{Conservation-laws-B}
\frac{\text{d}}{\text{d}t}\int_{\Omega_{(t)}} B \text{d}\Omega= \int_{\partial\Omega_{(t)}} B \textbf{V}^m \cdot \textbf{n} \text{d}s.
\end{equation}
The equation is derived based on the Reynolds transport theorem, similar to the construction of Eq. (\ref{Reynolds-transport-theorem}).
For the bottom topography function $B$, its local time partial derivative is zero.
The spatial discretization and one-stage second-order temporal discretization for Eq. (\ref{Conservation-laws-B}) is given as
\begin{equation}\label{Conservation-laws-B-discrete}
\begin{split}
B^{n+1}|\Omega^{n+1}|= B^{n}|\Omega^{n}| &+ \sum_{l,k}\big( \omega_k B_{l,k} \textbf{V}^m_{l,k} \cdot \textbf{n}_l \big) \big|\Gamma_{l} \big|\Delta t \\
&+\sum_{l,k} \big(\omega_k \textbf{V}^m_{l,k}\cdot \nabla B_{l,k} \textbf{V}^m_{l,k} \cdot \textbf{n}_l \big)\big|\Gamma_{l} \big|\frac{\Delta t^2}{2} \\
&+\sum_{l,k} \omega_k B_{l,k}(\mathbf{V}_1\times \mathbf{V}_2)\cdot \mathbf{k} \frac{\Delta t^2}{2}.
\end{split}
\end{equation}

\subsection{Cell-averaged gradients update on a moving mesh}

In this study, a compact GKS will be constructed for the SWE on a moving mesh. In the compact GKS, in addition to the cell-averaged flow variables given by Eq. (\ref{SWE-moving-cell-time}), the gradients of the cell-averaged flow variables are also updated via the Gauss's theorem,
\begin{equation}\label{slope}
\nabla \mathbf{W}_j^{n+1} =\frac{1}{|\Omega_j^{n+1}|} \int_{\partial \Omega_j^{n+1}} {\bf W}^{n+1}(\mathbf{x}) {\bf n} \mathrm{d} s,
\end{equation}
where the normal vector $\bf n$ is based on the moved cell $\Omega_j^{n+1}$ at $t^{n+1}$. The spatial discretization of Eq. (\ref{slope}) is obtained as
\begin{align*}
\nabla \mathbf{W}_j^{n+1} =\frac{1}{|\Omega_j^{n+1}|} \sum_{l,k} \omega_k {\bf W}^{n+1}_{l,k} {\bf n}\big|\Gamma_l\big|.
\end{align*}
On a moving mesh, the interface value ${\bf W}^{n+1}_{l,k}$ is determined by Eq. (\ref{slope-W-evolution}).
The second-order accurate time update for ${\bf W}^{n+1}_{l,k}$ is given by
\begin{equation}\label{slope-W-evolution-discrete}
\mathbf{W}^{n+1}_{l,k}=\mathbf{W}_{l,k} + \frac{\partial \mathbf{W}_{l,k}}{\partial t}\Delta t + \textbf{V}^m_{l,k} \cdot \nabla\textbf{W}_{l,k}\Delta t,
\end{equation}
where $\mathbf{W}_{l,k}$, $(\mathbf{W}_t)_{l,k}$, and $\nabla\textbf{W}_{l,k}$ can be obtained via direct numerical modeling \cite{zhao2023direct,zhao2021-swe} based on the time-dependent evolution solution of the GKS, presented in the Appendix.

It is worth noting that, in general, $\mathbf{W}^{n+1}_{l,k}$ takes different values, $\mathbf{W}^{n+1,l}_{l,k}$ and $\mathbf{W}^{n+1,r}_{l,k}$, within the cells on the two sides of an interface. When the flow is sufficiently smooth, the difference between $\mathbf{W}^{n+1,l}_{l,k}$ and $\mathbf{W}^{n+1,r}_{l,k}$ is consistent with the spatiotemporal discretization accuracy.
Providing two time-dependent, discontinuous evolving solutions at the mesh interface is key for high-order compact GKS to achieve compact spatial discretization and strong robustness. For popular finite volume schemes using approximate Riemann solvers and finite difference schemes, it is difficult to obtain the evolving solution $\mathbf{W}^{n+1}_{l,k}$ at the interface (or at half points).

\section{GCL-preserving and well-balanced property}

This section provides a theoretical analysis of the GCL-preserving and the well-balanced properties of the proposed ALE compact GKS on moving unstructured meshes.

\subsection{GCL-preserving property}

Numerical schemes on moving meshes necessarily include contributions arising from mesh motion. The space-time ALE scheme developed here is faithfully based on the mesh kinematics, enabling a consistent derivation of motion-induced terms.
For implementation simplicity, however, certain contributions are selectively omitted, which necessitates an assessment of the resulting accuracy. Preserving GCL is a fundamental criterion for evaluating moving mesh schemes.
Specifically, the GCL requires that Eq. (\ref{SWE-moving-cell-time}) hold identically for any prescribed mesh motion when the flow is steady and spatially uniform (i.e., $U$, $V$, $h$, and $B$ are constant in this study). Under these conditions, Eq. (\ref{SWE-moving-cell-time}) reduces to a purely geometric update of the cell area. Consequently, verification of the GCL reduces to checking whether the interface fluxes constructed from $\mathbf{L}_j^n$ and $\mathrm{d}\mathbf{L}_j^n/\mathrm{d}t$ coincide with the exact face-swept measure over the time step.

Considering uniform, steady flow and with the geometric conditions $\sum_{l=1}^3n_{l,x}|\Gamma_l|=0$ and $\sum_{l=1}^3n_{l,y}|\Gamma_l|=0$, there are $\mathcal{L}_{1,j}^n=0$ and $\mathrm{d} \mathcal{L}_{1,j}^n/\mathrm{d} t=0$.
At the same time, it is evident that $\mathcal{L}_{3,j}^n=0$ and $\mathrm{d} \mathcal{L}_{3,j}^n/\mathrm{d} t=0$.
On the RHS of Eq. (\ref{SWE-moving-cell-time}), the flux associated with $\mathcal{L}_{2,j}^n$ and $\mathrm{d}\mathcal{L}_{2,j}^n/\mathrm{d}t$ are
\begin{equation}\label{GCL-L2-operator}
\begin{split}
Tot(\mathcal{L}_{2,j}^n) \equiv & \mathcal{L}_{2,j}^n \Delta t + \frac{\mathrm{d} \mathcal{L}_{2,j}^n}{\mathrm{d} t}\frac{\Delta t^2}{2} \\
=& \sum_{l,k}\big( \omega_k \textbf{F}^m_{l,k} \cdot \textbf{n}_l \big) \big|\Gamma_{l} \big| \Delta t \\
&+ \sum_{l,k}\big( \omega_k (\partial\textbf{W}_{l,k}/\partial t+\mathbf{V}^m_{l,k}\cdot \nabla\mathbf{W}_{l,k}) \textbf{V}^m_{l,k} \cdot \textbf{n}_l \big) \big|\Gamma_{l} \big| \frac{\Delta t^2}{2} \\
&+ \sum_{l,k}\big( \omega_k \textbf{W}_{l,k} (\mathbf{V}_1\times \mathbf{V}_2)\cdot \mathbf{k} \big) \frac{\Delta t^2}{2},
\end{split}
\end{equation}
where $\textbf{F}^m= \textbf{W} \textbf{V}^m$.
Since the flux computation on each edge is independent, for simplicity, the following calculation considers only the flux on a single edge as shown in Fig. \ref{0-cell-moving}. Thus, the flux on one edge is
\begin{equation}\label{GCL-L2-operator-1}
\begin{split}
Tot_l(\mathcal{L}^n_{2,j})&=  \sum _{k} \omega_k \textbf{V}^m_k \cdot \textbf{n}_l  \big|\Gamma_{l} \big| \Delta t +  \sum _{k} \omega_k (\mathbf{V}_1\times \mathbf{V}_2)\cdot \mathbf{k} \frac{\Delta t^2}{2}\\
&= \frac{1}{2}(U_1^m+U_2^m)|\Gamma_l|\Delta t + (\mathbf{V}_1\times \mathbf{V}_2)\cdot \mathbf{k} \frac{\Delta t^2}{2}.
\end{split}
\end{equation}
On the other hand, for a line segment moving at the velocity of $\mathbf{V}^m$ as shown in Fig. \ref{0-cell-moving}, the quadrilateral area swept over a time step $\Delta t$ can be directly calculated as
\begin{equation}\label{mesh-moving-area}
\begin{split}
\Delta S&=\frac{1}{2}(U_1^m+U_2^m)|\Gamma_l|\Delta t + (\mathbf{V}_1\times\mathbf{V}_2)\cdot \mathbf{k} \frac{\Delta t^2}{2} \\
&=Tot_l(\mathcal{L}_{2,j}^n).
\end{split}
\end{equation}
Therefore, the flux $Tot_l(\mathcal{L}^n_{2,j})$ in Eq. (\ref{GCL-L2-operator-1}) is exactly equal to the area $\Delta S$ swept by the cell interfaces in Eq. (\ref{mesh-moving-area}), indicating that the current ALE scheme rigorously preserves GCL.
It should be noted that the right-hand term in Eq. (\ref{GCL-L2-operator-1}) can take positive and negative values. The positive and negative values correspond to the increase and decrease in the area of the cell, respectively, due to the moving of the interface. The positions of $\mathbf{V}_1$, $\mathbf{V}_2$, and the normal vector $\mathbf{n}_l$ satisfy the right-hand rule, meaning that when moving from point $A$ to point $B$, the normal vector $\mathbf{n}_l$ points towards the right side of the hand.

\subsection{Well-balanced property}

For the SWE with non-flat bottom, a central challenge is to enforce an exact discrete balance between convective fluxes and topographic source terms so that steady states are preserved, such as the lake-at-rest state; this is the well-balanced property.
In an ALE scheme on a moving mesh, achieving well-balancedness requires that the discretized source term be exactly balanced by two components of the numerical flux: (i) the physical hydrodynamic flux and (ii) the geometric flux associated with the interface-swept area induced by mesh motion.

To assess this property, we consider a two-dimensional stationary equilibrium with initial data given by a constant free-surface elevation $h + B$ and zero velocity $\mathbf{V} = \mathbf{0}$; the bathymetry $B$ is taken to be piecewise linear. 
The hydrodynamic flux is evaluated as
\begin{equation}
\begin{split}
Tot(L^n_{1,j})&=L^n_{1,j}\Delta t +\frac{\mathrm{d}L^n_{1,j}}{\mathrm{d}t} \frac{\Delta t^2}{2} \\
&=-\sum_{l,k}\frac{1}{4} \big(0,Gh^2_{l,k}n_{l,x},Gh^2_{l,k}n_{l,y} \big)^{\mathrm{T}}\big| \Gamma_l\big| \Delta t -\sum_{l,k}\frac{1}{2}Gh_{l,k}\big(0,h_x,h_y \big)^{\mathrm{T}}\mathbf{V}^m_{l,k}\cdot \mathbf{n}_l \big| \Gamma_l\big| \frac{\Delta t^2}{2}\\
&=-G\big|\Omega_j^n \big|\big(0,h_{j,x}h_j,h_{j,y}h_j\big)^{\mathrm{T}}\Delta t -\sum_{l,k}\frac{1}{2}Gh_{l,k}\big(0,h_x,h_y \big)^{\mathrm{T}}\mathbf{V}^m_{l,k}\cdot \mathbf{n}_l \big| \Gamma_l\big| \frac{\Delta t^2}{2}\\
&=G\big|\Omega_j^n \big|\big(0,B_{j,x}h_j,B_{j,y}h_j\big)^{\mathrm{T}}\Delta t +\sum_{l,k}\frac{1}{2}Gh_{l,k}\big(0,B_{j,x},B_{j,y} \big)^{\mathrm{T}}\mathbf{V}^m_{l,k}\cdot \mathbf{n}_l \big| \Gamma_l\big| \frac{\Delta t^2}{2}.
\end{split}
\end{equation}
The numerical flux due to mesh motion is
\begin{equation}
\begin{split}
Tot(L^n_{2,j})=&L^n_{2,j}\Delta t +\frac{\mathrm{d}L^n_{2,j}}{\mathrm{d}t} \frac{\Delta t^2}{2} \\
=&\sum_{l,k}\big(\frac{1}{2}(h,0,0)^\mathrm{T} \mathbf{V}^m_{l,k}\cdot \mathbf{n}_l \big)\big| \Gamma_l\big|\Delta t \\
&+\sum_{l,k}\big(\frac{1}{2}(\mathbf{V}^m_{l,k}\cdot \nabla h,0,0)^\mathrm{T} \mathbf{V}^m_{l,k}\cdot \mathbf{n}_l \big)\big| \Gamma_l\big| \frac{\Delta t^2}{2} \\
&+\sum_{l,k} \frac{1}{2} (h,0,0)^\mathrm{T} (\mathbf{V}_1\times \mathbf{V}_2)\cdot \mathbf{k} \frac{\Delta t^2}{2}.
\end{split}
\end{equation}
The source term is
\begin{equation}
\begin{split}
Tot(L^n_{3,j})=&L^n_{3,j}\Delta t +\frac{\mathrm{d}L^n_{3,j}}{\mathrm{d}t} \frac{\Delta t^2}{2} \\
=&-G h_j^n \big|\Omega_j^n \big|(0,B_{j,x}^n,B_{j,y}^n)^\mathrm{T}\Delta t \\
&-G L_{2,j}^{h} (0,B_{j,x}^n,B_{j,y}^n)^\mathrm{T} \frac{\Delta t^2}{2},
\end{split}
\end{equation}
and
\begin{align*}
\begin{split}
L_{2,j}^{n,h}=\sum_{l,k}\big(\frac{1}{2}h \mathbf{V}^m_{l,k}\cdot \mathbf{n}_l \big)\big| \Gamma_l\big|.
\end{split}
\end{align*}
There is $Tot(L^n_{1,j})+Tot(L^n_{3,j})=0$, so these contributions cancel exactly. Consequently, the updated momentum becomes zero and the water at rest is preserved, i.e.,
\begin{align*}
U^{n+1}_j=V^{n+1}_j=0.
\end{align*}
It remains to verify that the water surface elevation remains constant. The update of water height is given as
\begin{align*}
\begin{split}
h_j^{n+1}\big|\Omega_j^{n+1}\big|=&h_j^{n}\big|\Omega_j^{n}\big| +\sum_{l,k}\big(\frac{1}{2}h \mathbf{V}^m_{l,k}\cdot \mathbf{n}_l \big)\big| \Gamma_l\big|\Delta t \\
&+\sum_{l,k}\big(\frac{1}{2}\mathbf{V}^m_{l,k}\cdot \nabla h \mathbf{V}^m_{l,k}\cdot \mathbf{n}_l \big)\big| \Gamma_l\big| \frac{\Delta t^2}{2} +\sum_{l,k} \frac{1}{2} h (\mathbf{V}_1\times \mathbf{V}_2)\cdot \mathbf{k} \frac{\Delta t^2}{2}
\end{split}
\end{align*}
Combining the above equation with Eq. (\ref{GCL-L2-operator-1}) for the bottom topography update, we obtain the update for the water surface elevation as
\begin{align*}
\begin{split}
(h_j^{n+1}+B_j^{n+1})\big|\Omega_j^{n+1}\big|=& h_{\mathrm{surf}}\big|\Omega_j^{n}\big| +h_{\mathrm{surf}}\sum_{l,k}\frac{1}{2}\mathbf{V}^m_{l,k}\cdot \mathbf{n}_l \big| \Gamma_l\big|\Delta t +h_{\mathrm{surf}}\sum_{l,k}\frac{1}{2}(\mathbf{V}_1\times \mathbf{V}_2)\cdot \mathbf{k} \frac{\Delta t^2}{2} \\
=&h_{\mathrm{surf}}\big|\Omega_j^{n}\big| +h_{\mathrm{surf}} \sum_{l,k}\frac{1}{2}\Delta S_{l} \\
=&h_{\mathrm{surf}}\big|\Omega_j^{n+1}\big|.
\end{split}
\end{align*}
Therefore, the constant water surface is preserved, i.e., $(h_j^{n+1}+B_j^{n+1})=h_{\mathrm{surf}}$.

\section{High-order compact reconstruction}

In this study, a fourth-order compact spatial reconstruction is employed to determine the values of flow variables at the cell interface, which are then used to compute the evolution solutions for these flow variables and their fluxes. To handle discontinuous solutions that may arise in the flow, a nonlinear reconstruction is necessary, for which a linear high-order reconstruction serves as a prerequisite. This study will adopt the previously developed fourth-order compact reconstruction based on the WENO method \cite{zhaocompact_tri}, and the corresponding formulas will be briefly presented in this section.

\subsection{Linear high-order compact reconstruction}

\begin{figure}[!htb]
\centering
\includegraphics[width=0.35\textwidth]{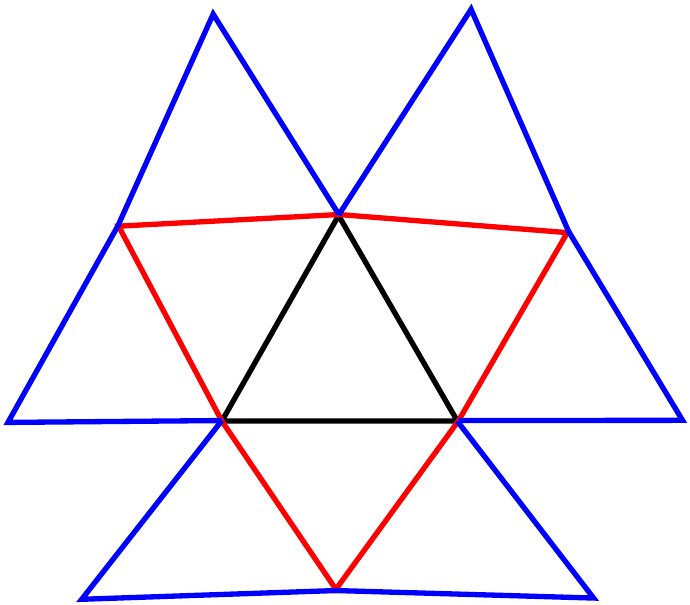}
\caption{\label{0-stencil-2d} A schematic of reconstruction stencil of compact GKS. The flow variables and their gradients are known in each cell.}
\end{figure}

The numerical scheme presented in Section 2 of this study updates the cell-averaged values of the flow variables and their gradients. Using these values, a fourth-order spatial reconstruction can be achieved on a compact stencil that involves only neighboring cells. Fig. \ref{0-stencil-2d} shows the compact stencil on a triangular mesh for implementing this fourth-order reconstruction.
The reconstruction cell, outlined in black, is given the local index $0$. Its three edge-neighbors and six vertex-neighbors are indexed sequentially from $1$ to $9$. In some special cases, a few of the vertex-neighboring cells may not exist.

For the fourth-order reconstruction, a cubic polynomial $P^3(\bm{x})$ is constructed on each cell as
\begin{align*}
\begin{split}
P^3(\bm{x})=\sum_{k=1}^{10} a_k \varphi_k(\bm{x}),
\end{split}
\end{align*}
where $a_k$ is the degrees of freedom (DOFs) of $P^3(\bm{x})$, $\varphi_k(\bm{x})$ denotes the complete set of cubic polynomial basis functions.
The coefficients $a_k$ of $P^3$ are calculated based on the condition that the cell averages of $P^3$ and the gradients over the stencil cells satisfy the given values, as follows:
\begin{align}\label{Recons-eqs}
\begin{split}
&\big(\frac{1}{\big|\Omega_i \big|} \int_{ \Omega_i} \varphi_k(\bm{x}) \mathrm{d}x \mathrm{d}y \big) a_k=Q_i, ~i=0,1,\cdots,9, \\
&\big(\frac{\Delta X}{\big|\Omega_i \big|} \int_{ \Omega_i} \varphi_{k,l}(\bm{x}) \mathrm{d}x \mathrm{d}y \big) a_k=Q_{i,l} \Delta X, ~i=0,1,2,3,~l=x,y, \\
\end{split}
\end{align}
where the same subscript $k$ of $\varphi_k$ and $a_k$ on the left-hand side of the equations follow the Einstein summation, $Q$ is any component of $\mathbf{W}$, and $\Delta X$ is the cell size (e.g., the square root of the cell area).
The second equation in Eq. (\ref{Recons-eqs}) is multiplied by $\Delta X$ to scale the terms on its left-hand side to the same order of magnitude as those in the first equation.

Considering the inconsistency between the number of coefficients $a_k$ and the number of reconstruction conditions in Eq. (\ref{Recons-eqs}), and in order to satisfy the conservation of $P^3$ on the reconstruction cell, a constrained least-squares method is employed. The resulting linear system of equations for $a_k$ is given as:
\begin{equation} \label{CLS-system}
\left(
\begin{array}{cc}
\mathbf{A}_{0,m} & 0 \\
2\mathbf{A}_{i,k} \mathbf{A}_{i,m} & \mathbf{A}_{0,m}^\mathrm{T} \\
\end{array}
\right)
\left(
\begin{array}{c}
\mathbf{a}_k \\
b \\
\end{array}
\right)
=
\left(
\begin{array}{cc}
1 & \mathbf{0} \\
\mathbf{0} & 2\mathbf{A}_{i,k} \\
\end{array}
\right)
\mathbf{q}.
\end{equation}
where $i=1,2,\cdots,17$.
The matrix $\mathbf{A}$ is formed from the cell averages of the basis functions and their derivatives from Eq. (\ref{Recons-eqs}). $\mathbf{A}$ is an $18\times10$ matrix, where its first $10$ rows are derived from the first equations in Eq. (\ref{Recons-eqs}), and the remaining $8$ rows are derived from the second. The vector $\mathbf{q}$ is an $18\times1$ column vector whose first ten components are $Q_i$ and whose last eight components are $Q_{i,l} \Delta X$ $(i=0,1,2,3,~ l=x,y)$. The variable $b$ is an auxiliary Lagrange multiplier introduced for the constraint.

\subsection{Nonlinear reconstruction }

To handle discontinuous solutions, a non-linear reconstruction procedure is required. The present study adopts the WENO method, which has been previously developed for the compact GKS. The fundamental principle of the WENO method is to adaptively combine several lower-order polynomials which are reconstructed on the sub-stencils.
In the compact GKS, the second-order reconstruction is defined on sub-stencils that each involve only one neighboring cell, and three such sub-stencils are possible. The linear polynomial $P^1_m$ is given by
\begin{align*}
\begin{split}
P^1_m(\bm{x})=Q_0 +a_{m,1}(x-x_0) +a_{m,2}(y-y_0), ~m=1,2,3,
\end{split}
\end{align*}
and $P^1_m$ satisfies the following conditions on the sub-stencil
\begin{align}\label{Recons-eqs-p1}
\begin{split}
&\big(\frac{1}{\big|\Omega_i \big|} \int_{ \Omega_i} \varphi_k(\bm{x}) \mathrm{d}x \mathrm{d}y \big) c_k=Q_i, ~i=m, \\
&\big(\frac{\Delta X}{\big|\Omega_i \big|} \int_{ \Omega_i} \varphi_{k,l}(\bm{x}) \mathrm{d}x \mathrm{d}y \big) c_k=Q_{i,l} \Delta X, ~i=m,~l=x,y, \\
\end{split}
\end{align}
By applying a constrained least-squares method with the first equation in Eq. (\ref{Recons-eqs-p1}) as the constraint, a linear system for the coefficients $c_k$ can be constructed. This system is analogous to Eq. (\ref{CLS-system}), and the details is therefore omitted for brevity.

Finally, the WENO reconstruction yields
\begin{equation}\label{4th-compact-HLP}
\begin{split}
&R(x)= \sum_{k=1}^{3}w_k P^1_k(\mathbf{x}) + w_0 \big( \frac{1+C}{C}P^3(\mathbf{x}) -\sum_{k=1}^{3} \frac{C_k}{C}P^1_k(\mathbf{x}) \big),
\end{split}
\end{equation}
where $w_k$ ($k=0,1,2,3$) are the nonlinear weights, and $C$ and $C_k$ are constant parameters set to $5$ and $1/3$, respectively. The details of the nonlinear reconstruction procedure can be found in \cite{zhao2023direct}.

\section{Mesh moving method}
In this study, we propose a simple algorithm for mesh adaptive moving.
The basic idea of this method is consistent with some existing adaptive mesh refinement methods, which involve placing more mesh cells in regions with larger gradients or second derivatives of the flow variables.
The innovation of this study lies in providing a simpler implementation approach, which we refer to as the variation-antigradient diffusion method.
Its key feature is that the update of each mesh node's position is performed independently, without the need for information on local mesh connectivity.
The algorithm is explicit and does not require matrix solving.

\begin{figure}[!htb]
\centering
\includegraphics[width=0.35\textwidth]{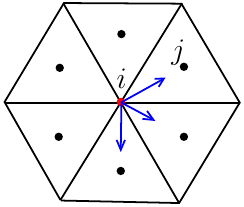}
\caption{\label{0-moving-adaptive-method} A schematic of mesh node motion. Adaptive mesh moving is achieved through node motion.}
\end{figure}

The algorithm for adaptive mesh node motion is outlined below. A schematic is provided in Fig. \ref{0-moving-adaptive-method}.
To balance computational efficiency, mesh motion is performed only every $N_m$ steps. In all numerical examples in this study, we set $N_m=4$.
\begin{enumerate}
    \item For each cell, the variation of the flow variable is quantified by the smoothness indicator of the cubic polynomial computed during the reconstruction step, denoted as $Var_j=IS(P^3)$.
    \item For each node $i$, the nodal variation $\overline{Var}_i$ is computed by averaging the variations from its $J$ neighboring cells, i.e., $\overline{Var}_i=\sum_{j=1}^J Var_j/J$.
    \item Calculate the potential displacements of node $i$ towards the centroids of its neighboring cells. For each neighboring cell $j$, the displacement is defined as $\Delta \mathbf{S}_{ij}=C \overrightarrow{P_iP_j} \Delta Var_{ij}$, where $\Delta Var_{ij}=\mathrm{Max}\{\epsilon,Var_{j}-\overline{Var}_{i}\}$, and $C$ is a parameter, set to $C=1/4$. The small parameter $\epsilon$ is assigned a value of $10^{-15}$.
    \item Determine the final displacement of node $i$ as $\Delta \mathbf{S}_i=\sum_{j=1}^J \Delta \mathbf{S}_{ij}/\sum_{j=1}^J \Delta Var_{ij}$, where the sums are taken over all $J$ neighboring cells. $\Delta \mathbf{S}_i$ is further limited by the local minimum cell size: if $\Delta \mathbf{S}_i<\Delta X_{min}$, then $\Delta \mathbf{S}_i=0$, where $\Delta X_{min}$ is the minimum size of cell $i$ and its neighbors $j$. The cell size is defined by the inradius.
    \item Mesh smoothing. After $N_s$ steps of adaptive mesh motion, the mesh is smoothed to prevent mesh singularities. In the smoothing step, the displacement of a mesh node is given by $\Delta \mathbf{S}_i=\sum_{j=1}^{J}\mathbf{X}_j-\mathbf{X}_i$. In the numerical experiments, $N_s$ is uniformly set to $6$.
\end{enumerate}

Additionally, to verify the scheme, a prescribed position- and time-dependent mesh motion is also employed, and the details are provided in the examples. Unlike the adaptive mesh motion, under the prescribed motion the mesh is moved and updated at every time step during the computation.

\section{Numerical validations}

In this section, the high-order compact GKS on moving meshes is applied to a series of numerical test cases to evaluate its ability to preserve the GCL and well-balanced property, and to capture discontinuities without spurious oscillations.
In the computation, the time step used is determined by the CFL condition with CFL$=0.5$.
The gravitational acceleration is taken as $G=1.0$ if not specified.
The collision time $\tau$ in the gas distribution function for inviscid flow at a cell interface is defined by
\begin{align*}
\tau= C_{num}\big|\frac{h_l^2-h_r^2}{h_l^2+h_r^2}\big|\Delta t,
\end{align*}
where $C_{num}=10$, and $h_l$ and $h_r$ are the water heights on the left and right sides of the cell interface.
The squared values of $h_l$ and $h_r$ represent the hydrostatic pressure. The reason for including the pressure jump term in the relaxation time is to enhance the artificial dissipation in case of bore waves.

\begin{figure}[!htb]
\centering
\includegraphics[width=0.475\textwidth]{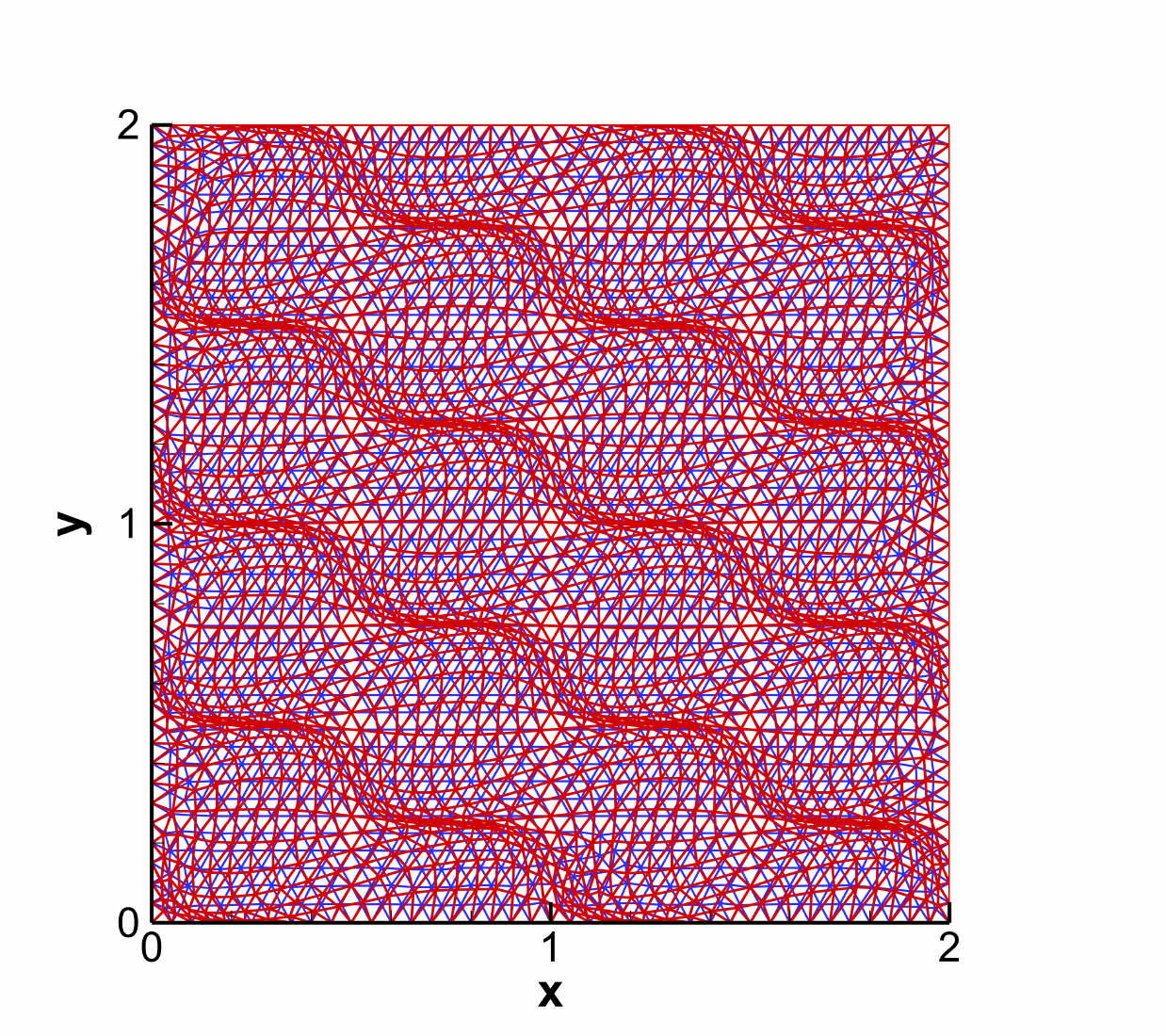}
\includegraphics[width=0.475\textwidth]{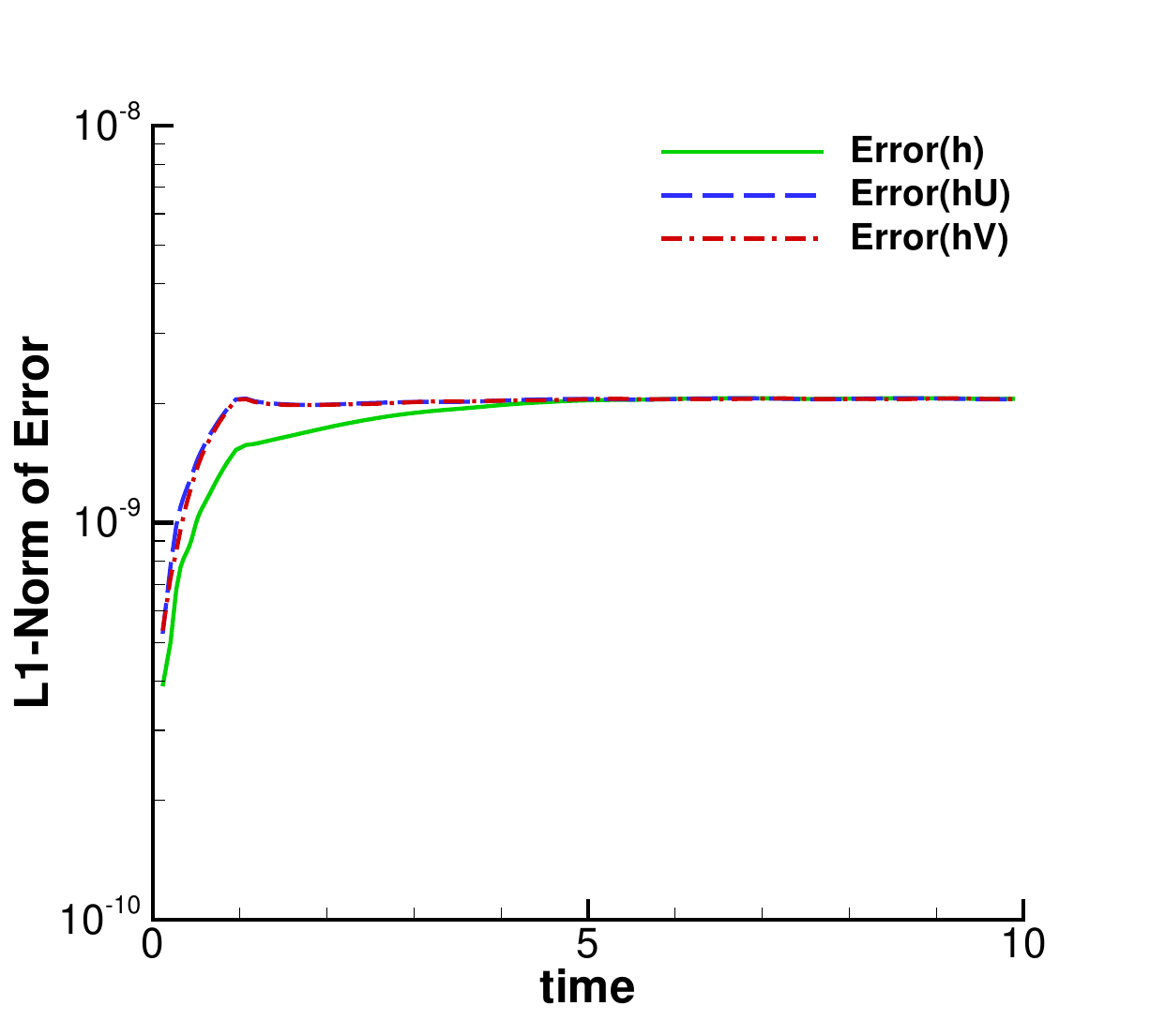}
\caption{\label{2d-geometric-law} Verification of the GCL of the compact GKS. Left: Comparison of the moving mesh at its initial state ($t=0$, blue) and a deformed state ($t=5.5$, red). Right: Time evolution of the numerical errors for water height and momentum.}
\end{figure}

\subsection{Verification of the GCL for the compact GKS}
The uniform free-stream problem is commonly used to verify the GCL-preserving property of a numerical scheme on moving meshes.
The initial condition of the uniform free-stream problem is set as
\begin{align*}
&h=1, \\
&(U,V)=(1,-1).
\end{align*}
The bottom is flat.
The computational domain is taken as $[0,2]\times[0,2]$, and the initially uniform triangular mesh with cell size $\Delta X=0.05$ is used.
The nodal coordinates of the moving mesh are defined as
\begin{align*}
\mathbf{x}=\mathbf{x}_0 +0.075\mathrm{sin}(\pi t)\mathrm{sin}(2\pi x)\mathrm{sin}(4\pi y)\mathbf{e},
\end{align*}
where $\mathbf{e}=(1,1)$ and $\mathbf{x}_0$ is the coordinates of nodes at $t=0$.

The left panel of Fig. \ref{2d-geometric-law} illustrates the significant deformation and displacement of the mesh at $t=5.5$ relative to its initial configuration.
The error histories for water height and momentum are plotted in the right panel, showing that the errors remain at a small, stable magnitude without growth over time.
This test case therefore confirms that the proposed compact GKS on moving meshes preserves the GCL on arbitrarily moving triangular meshes.

\begin{figure}[!htb]
\centering
\includegraphics[width=0.475\textwidth]{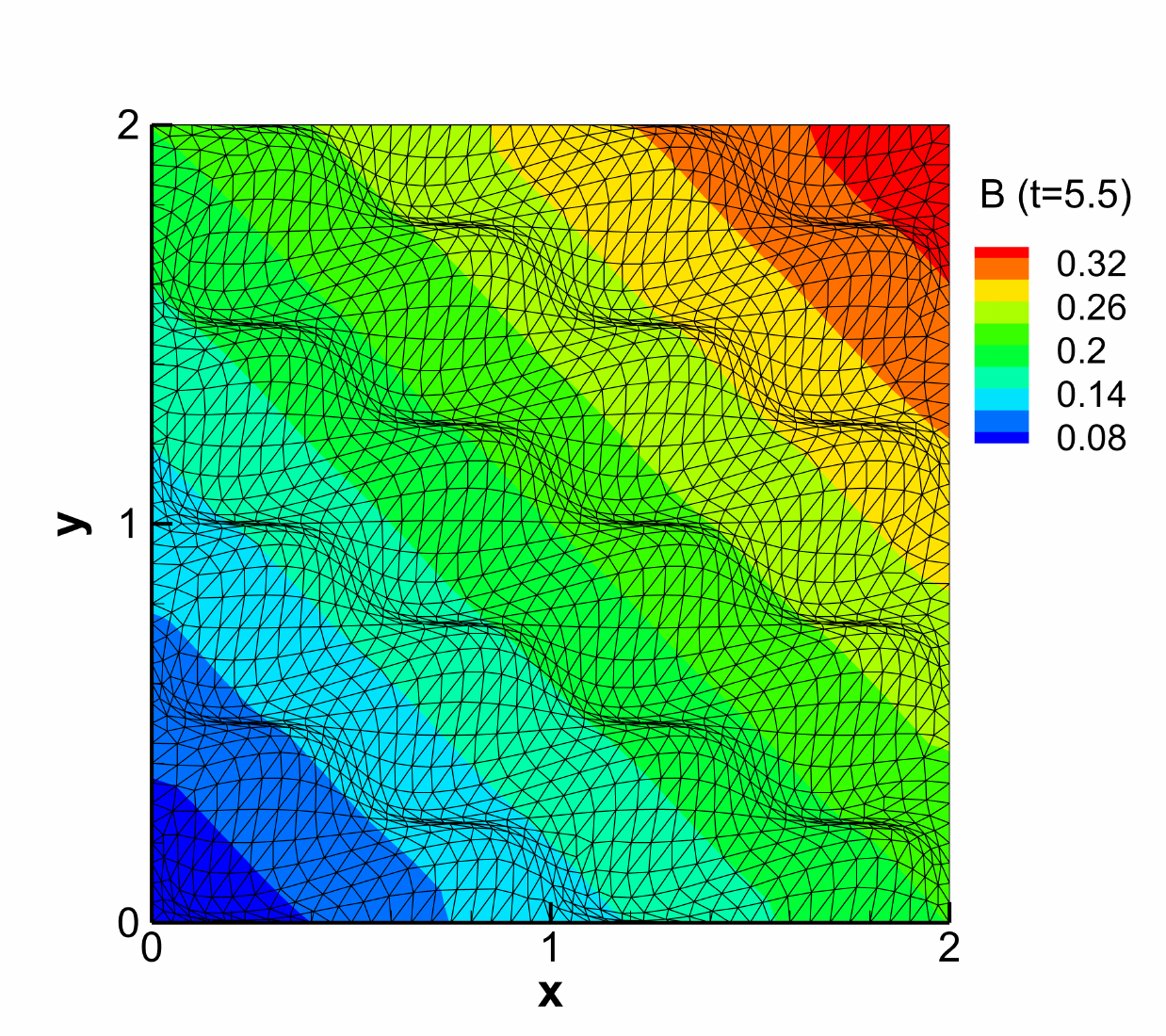}
\includegraphics[width=0.475\textwidth]{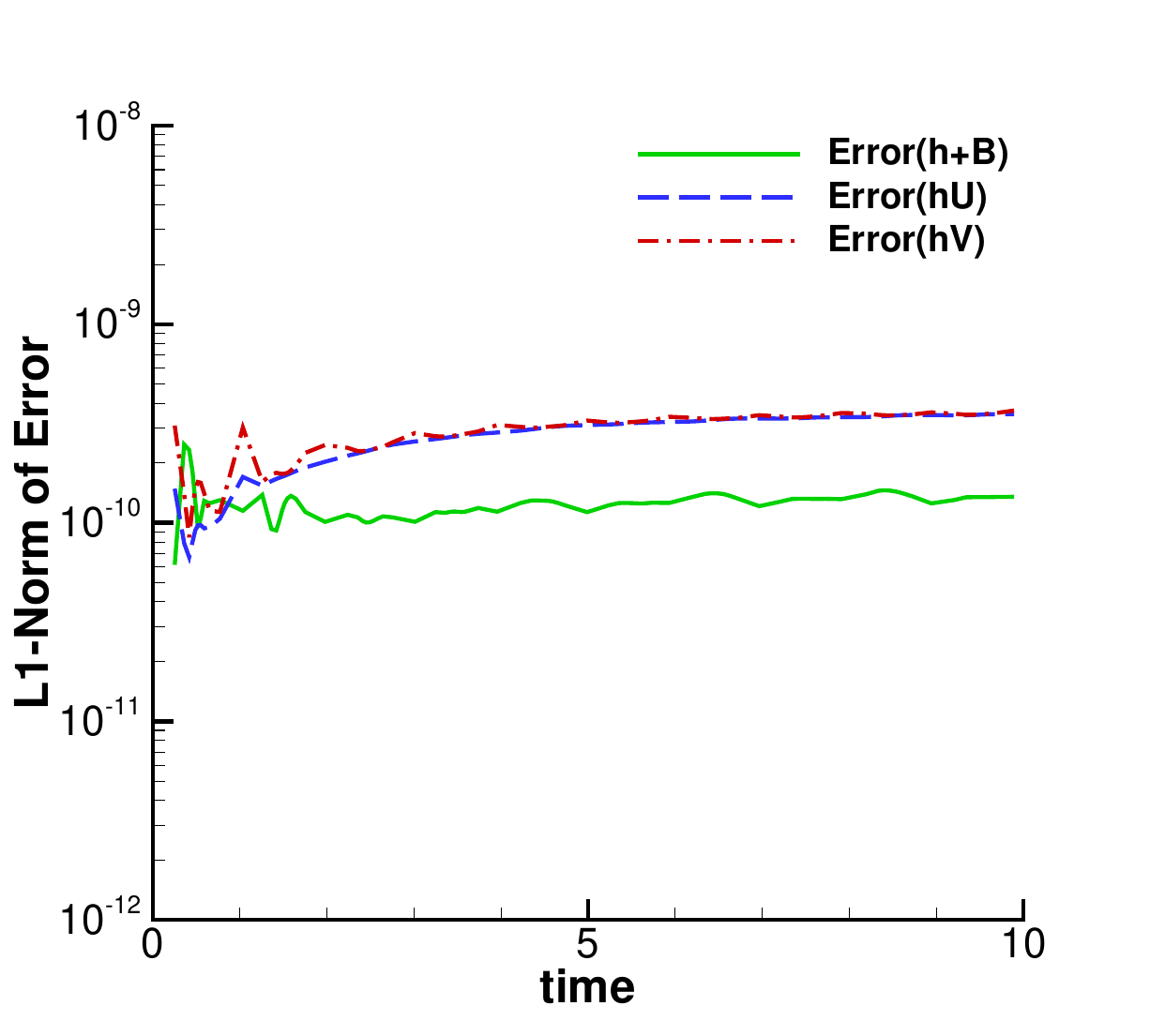}
\caption{\label{2d-balance-B-h-linear} Numerical verification of well-balanced property with a linear bottom topography. Left: Mesh and bottom profile at the time ($t=5.5$) of maximum mesh deformation. Right: Error history of water height and momentum.}
\end{figure}

\begin{figure}[!htb]
\centering
\includegraphics[width=0.475\textwidth]{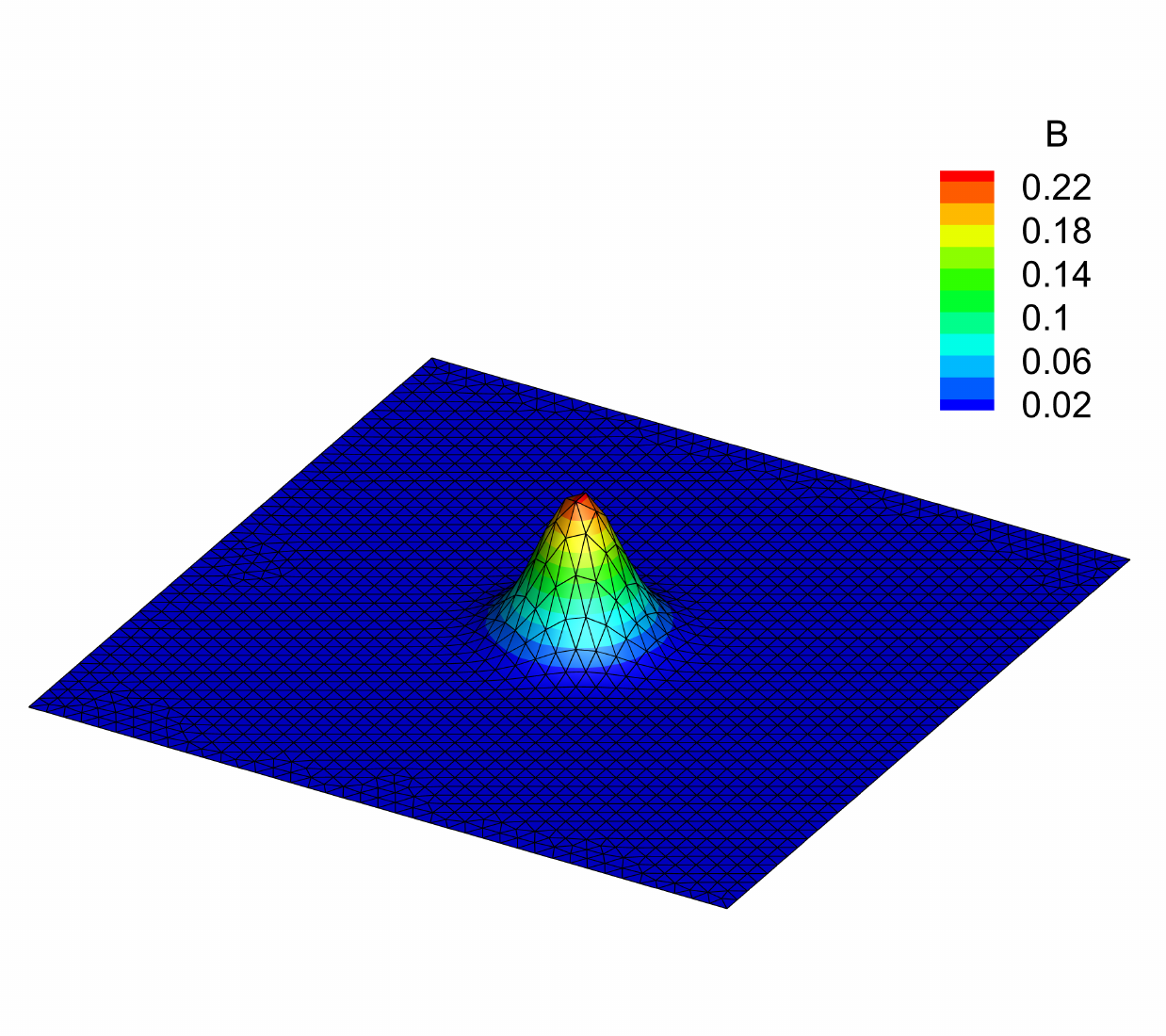}
\includegraphics[width=0.475\textwidth]{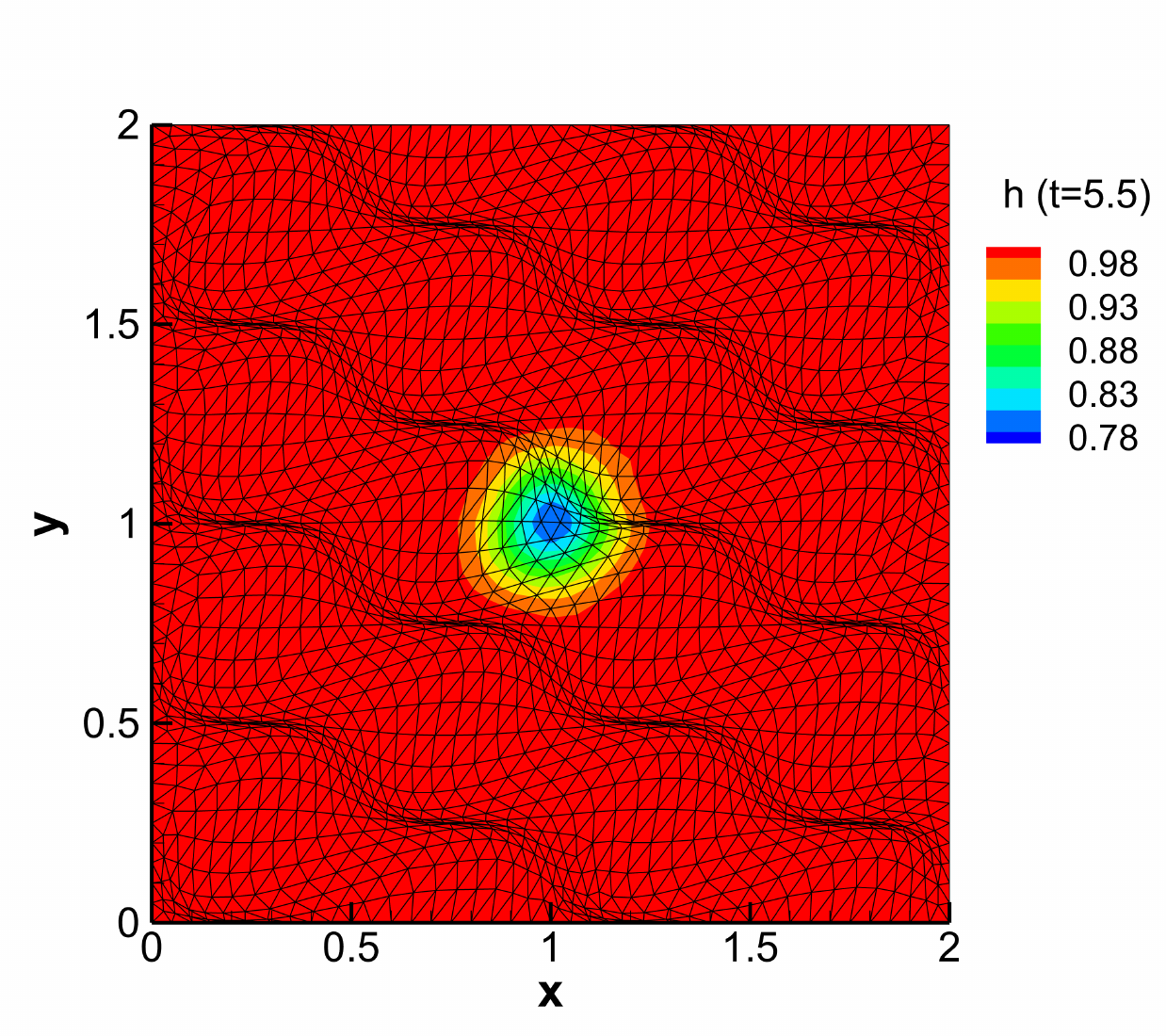}
\caption{\label{2d-balance-B-h} Numerical verification of well-balanced property with nonlinear bottom topography. Left: The discretized bottom profile at $t=0$. Right: Distributions of mesh and water height at $t=5.5$.}
\end{figure}

\begin{figure}[!htb]
\centering
\includegraphics[width=0.475\textwidth]{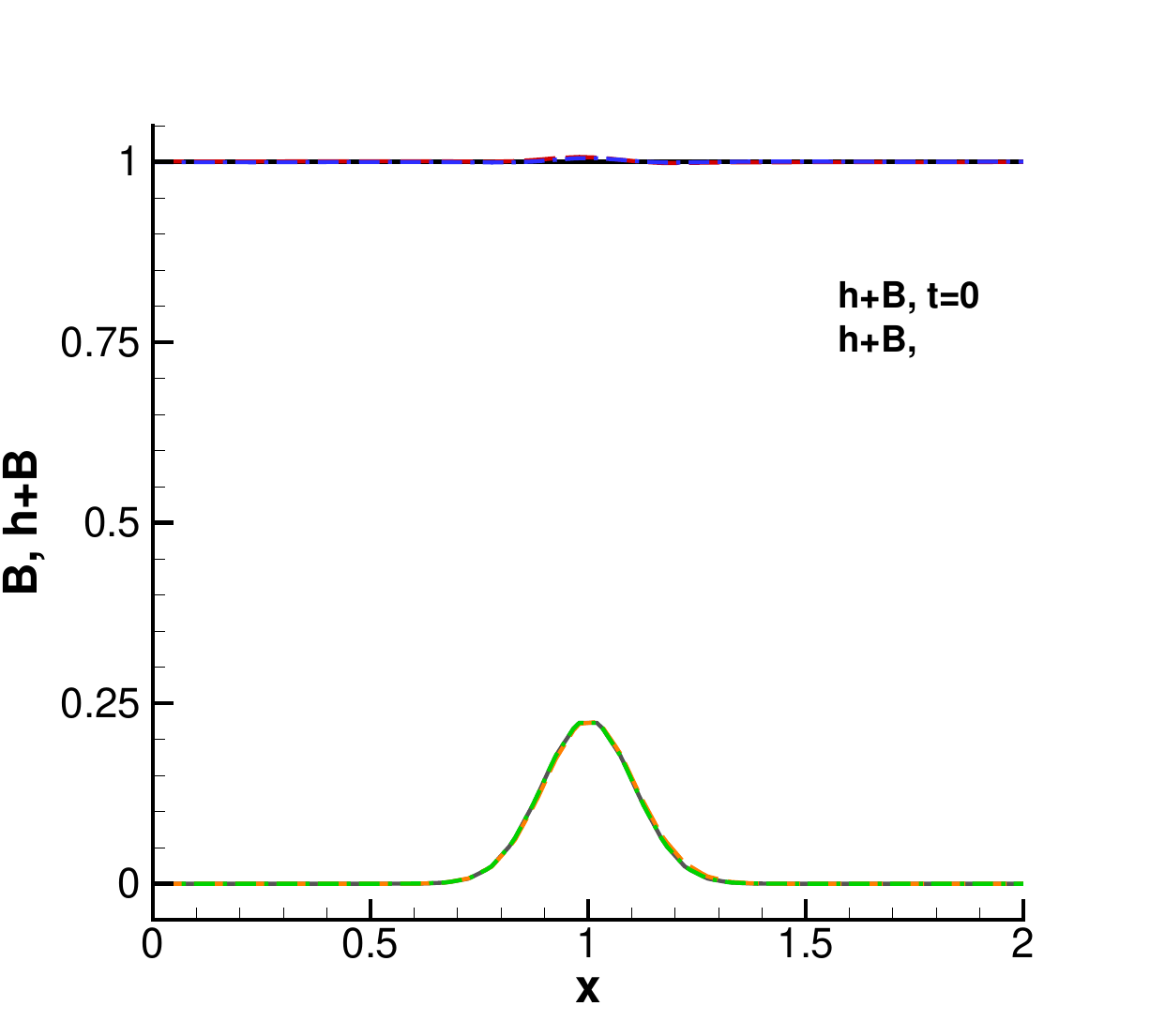}
\includegraphics[width=0.475\textwidth]{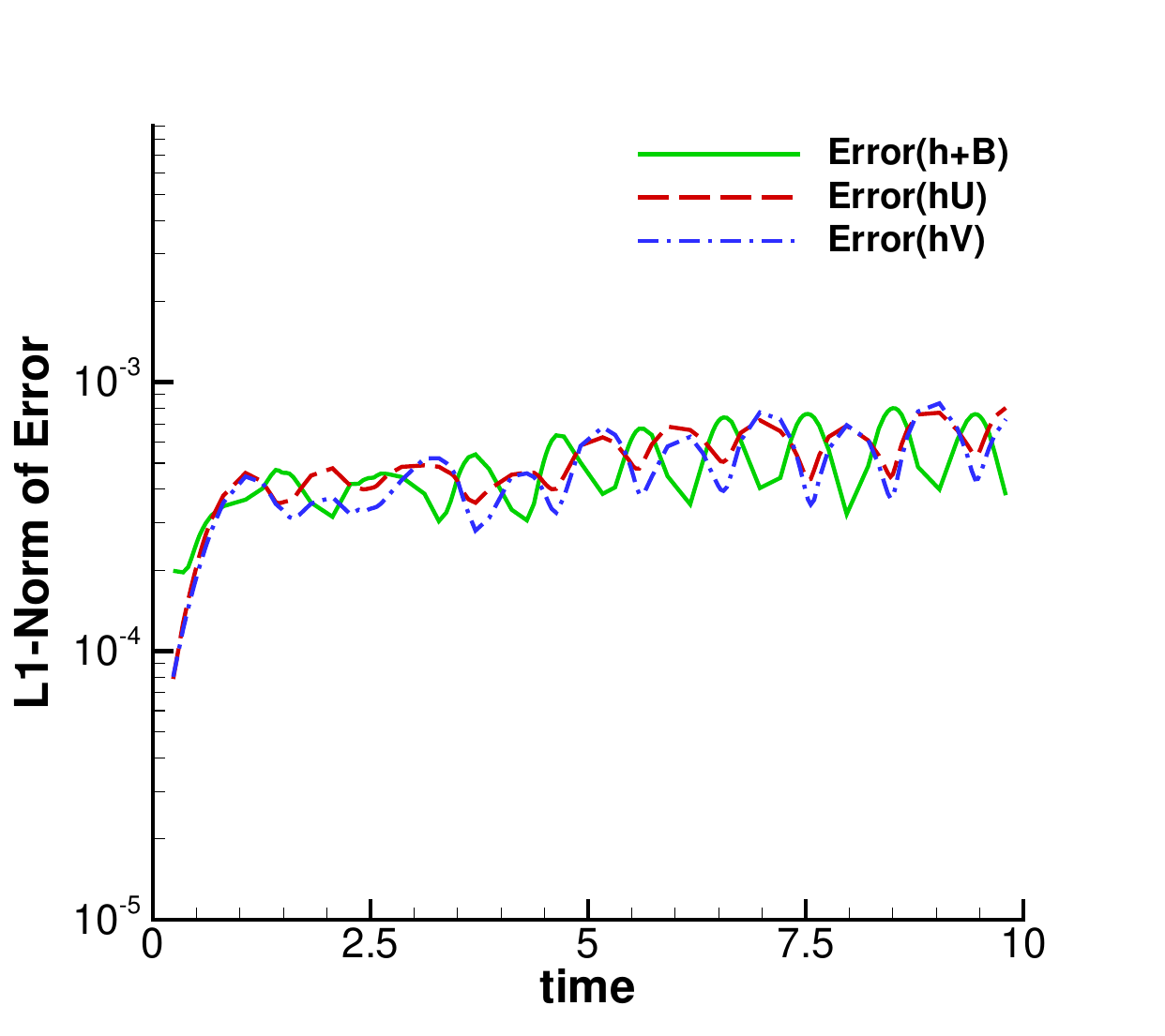}
\caption{\label{2d-balance-error} Numerical verification of well-balanced property with nonlinear bottom topography. Left: Distributions of water surface and bottom topography along $y=1$ at different times. Right: The time history of numerical errors.}
\end{figure}

\subsection{Verification of the well-balanced property for the compact GKS}
The well-balanced property of the compact GKS on moving meshes under non-flat bottom is verified using a lake-at-rest steady-state problem. The solution of the steady-state problem is given by
\begin{align*}
&h=1-B(x,y), \\
&(U,V)=(0,0).
\end{align*}
The computations are performed on the domain $[0,2]\times[0,2]$, with boundary conditions prescribed by the given lake-at-rest steady-state solution.
Two types of bottom topography are considered. The first is a linear bottom distribution given by
\begin{align*}
B(x,y)=0.05+0.075(x+y).
\end{align*}
The second is a more complex nonlinear distribution defined by
\begin{align*}
B(x,y)=0.25e^{-50[(x-1)^2+(y-1)^2]}.
\end{align*}
The distinction between these two bottom topographies in terms of their impact on the numerical scheme lies in whether, after updating the cell-averaged bottom values, one can exactly reconstruct a piecewise linear bottom distribution within each mesh cell.
The computational mesh and mesh motion are identical to those used in the GCL verification case. The mesh vertices evolve in time as
\begin{align*}
\mathbf{x}=\mathbf{x}_0 +0.075\mathrm{sin}(\pi t)\mathrm{sin}(2\pi x)\mathrm{sin}(4\pi y)\mathbf{e},
\end{align*}
where $\mathbf{e}=(1,1)$ and $\mathbf{x}_0$ is the coordinates of nodes at $t=0$.

In the linear bottom topography case, Fig. \ref{2d-balance-B-h-linear} (left) illustrates the mesh and bottom profile at $t=5.5$, the time of maximum mesh deformation. The error histories for water surface height and momentum are shown in the right panel, stabilizing below $10^{-9}$. This numerically validates the well-balanced property of the compact GKS on moving meshes.

For the nonlinear bottom case, the initial discretized bottom and the water surface at $t=5.5$ are shown in Fig. \ref{2d-balance-B-h}. Fig. \ref{2d-balance-error} (left) provides a quantitative comparison of the water surface and bottom topography along $y=1$ at various times, indicating the bottom is accurately preserved while the maximum relative error in water surface height remains within $1\%$. The error histories in Fig. \ref{2d-balance-error} (right) show that errors in water surface height and momentum stabilize below $10^{-3}$.
This larger error for the nonlinear bottom case arises because the reconstruction cannot maintain a continuous piecewise linear distribution for the bottom at cell interfaces after it is updated on the moving cells.

\subsection{2-D propagation of perturbation in a steady state}
This test case is about capturing the two-dimensional small perturbation in steady state.
The bottom topography is set as
\begin{align*}
& B(x,y)=0.8e^{-5(x-0.9)^2-50(y-0.5)^2}.
\end{align*}
The initial steady state with small perturbation is
\begin{equation*}
h(x,y) = \begin{cases}
1-B(x,y) +0.01, ~~ & 0.05 \leq x \leq 0.15,\\
1-B(x,y),  & \mathrm{otherwise}.
\end{cases}
\end{equation*}
and the velocity is $U=V=0$.
The computational domain takes $[0,2]\times[0,1]$.
Free boundary condition is applied to the left and right boundaries, while slip wall boundary condition is used for the top and bottom boundaries.
The triangular mesh with the initial cell size $\Delta X=1/100$ is used.
The nodal coordinates of the moving mesh are given as
\begin{align*}
\mathbf{x}=\mathbf{x}_0 +0.04\mathrm{sin}(\pi t)\mathrm{sin}(2\pi x)\mathrm{sin}(4\pi y)\mathbf{e}.
\end{align*}

\begin{figure}[!htb]
\centering
\includegraphics[width=0.495\textwidth]{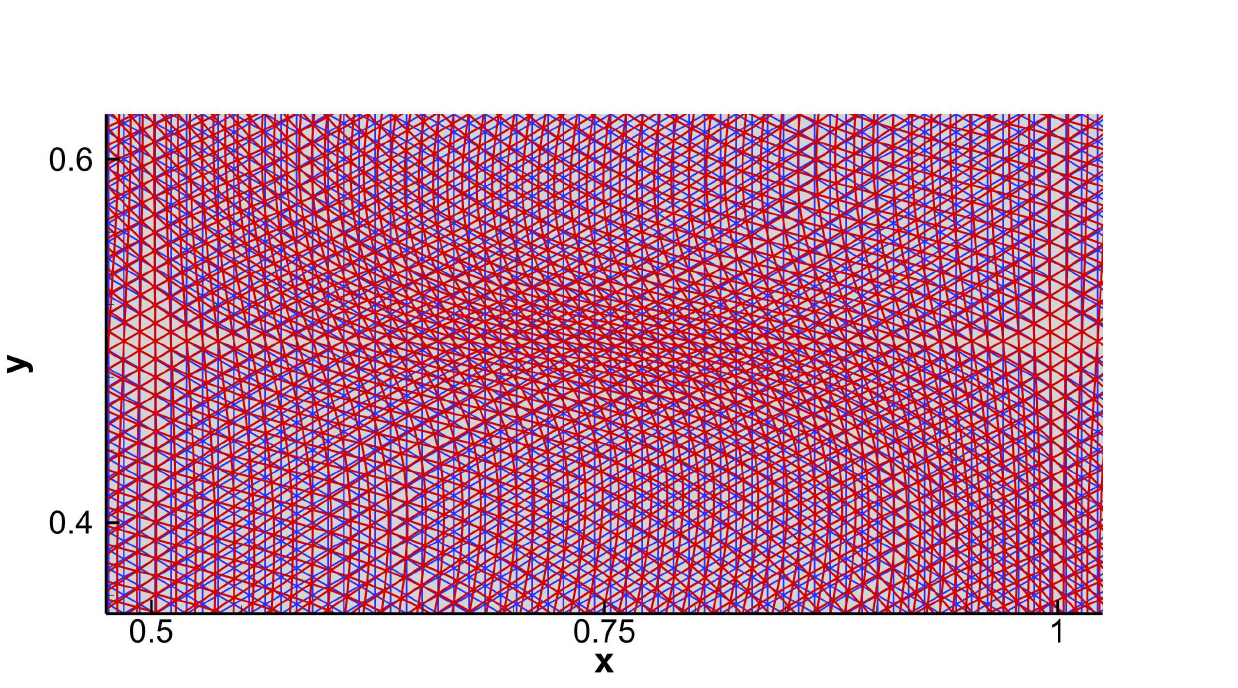}
\includegraphics[width=0.495\textwidth]{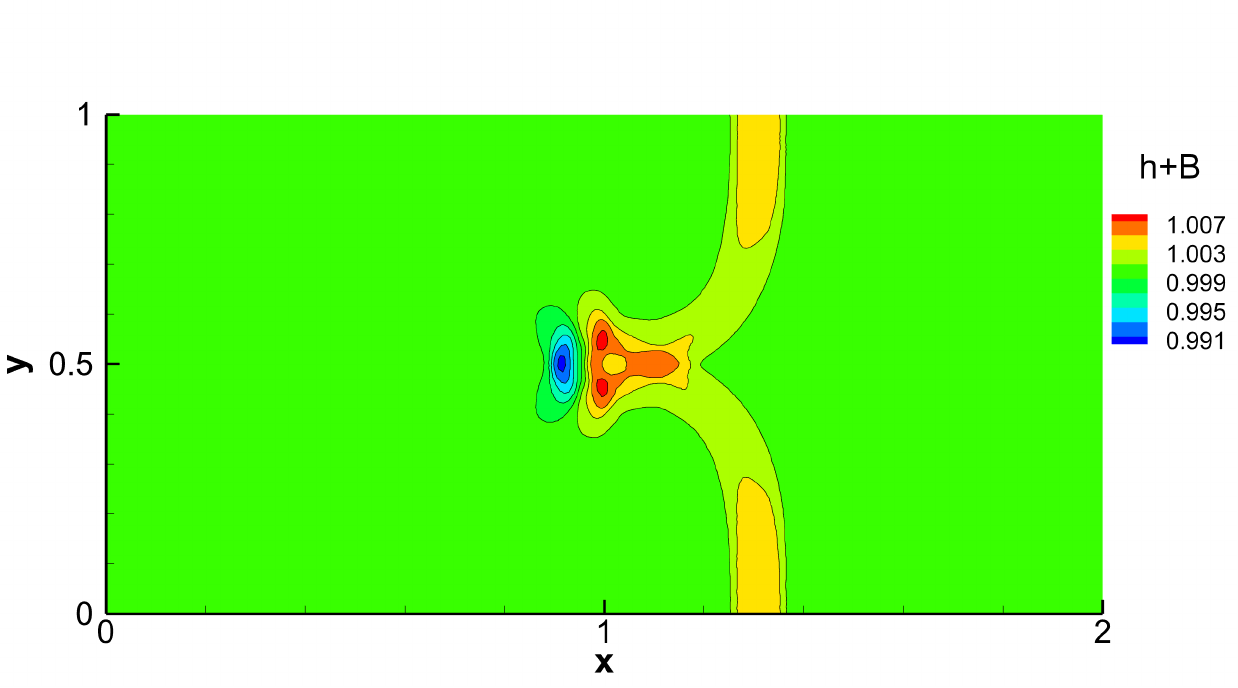}
\caption{\label{2d-nonflat-perturbation-1} 2-D propagation of a perturbation over non-flat bottom topography. Left: Comparison of the computational mesh at its initial state ($t=0$, shown in blue) and its deformed state at $t=0.5$ (shown in red). Right: Contours of the final water surface elevation at $t=1.2$.}
\end{figure}

\begin{figure}[!htb]
\centering
\includegraphics[width=0.495\textwidth]{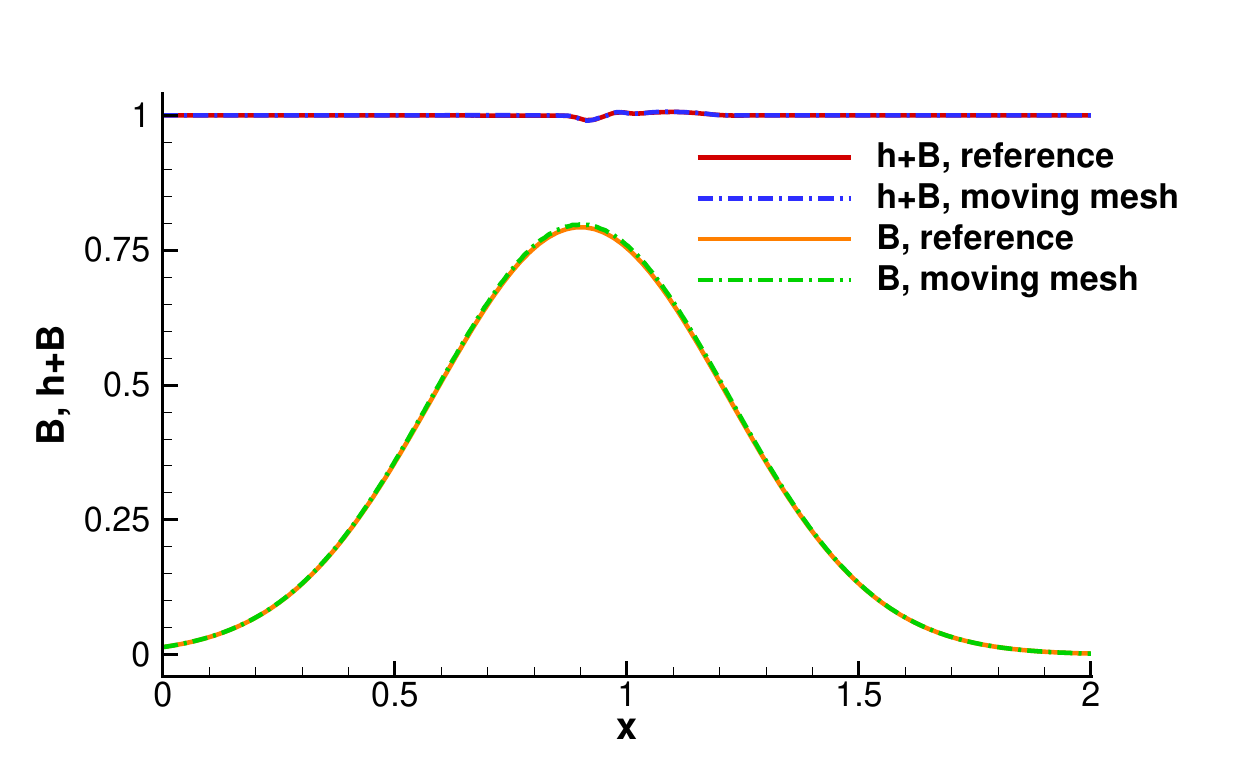}
\includegraphics[width=0.495\textwidth]{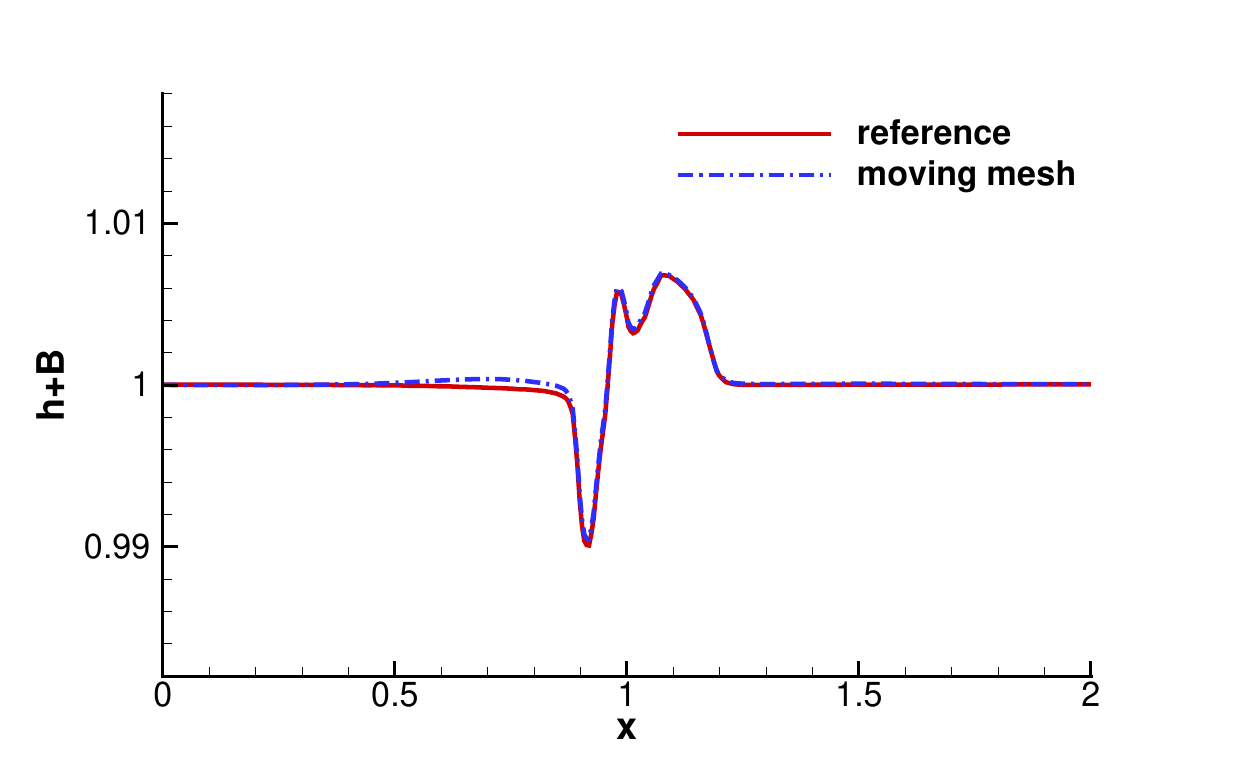}
\caption{\label{2d-nonflat-perturbation-2} 2-D propagation of a perturbation over non-flat bottom topography. Left:  Profiles of the water surface elevation and the underlying bottom topography along the horizontal centerline at $t=1.2$. Right: A zoomed-in view of the small perturbations in the water surface elevation.}
\end{figure}

Fig. \ref{2d-nonflat-perturbation-1} displays the computational mesh at $t=0.5$ and the contours of water surface elevation at $t=1.2$.
At $t=0.5$, the mesh nodes have the maximum displacement relative to their initial positions.
The numerical results demonstrate that small perturbations in the water surface elevation are accurately captured on the moving mesh, which validates the accuracy of the present scheme in handling the combined effects of bottom topography and mesh motion.
Fig. \ref{2d-nonflat-perturbation-2} provides a quantitative comparison of the bottom topography and the water surface elevation at the end of the computation against a reference solution.
The reference solution is obtained using a fourth-order compact GKS on the fixed initial mesh of this test case.
As observed in the left panel of Fig. \ref{2d-nonflat-perturbation-2}, a small discrepancy ($<6\times 10^{-3}$) exists in the computed bottom profile near the crest of the curve, which is attributed to the nonlinear distribution of the topography. The right panel of Fig. \ref{2d-nonflat-perturbation-2} presents a zoomed-in view of the small perturbations in the water surface elevation. Apart from a minor error ($< 6\times 10^{-4}$) in the vicinity of the bottom crest, the numerical solution shows excellent agreement with the reference solution across the rest of the domain.

\begin{figure}[!htb]
\centering
\includegraphics[width=0.475\textwidth]{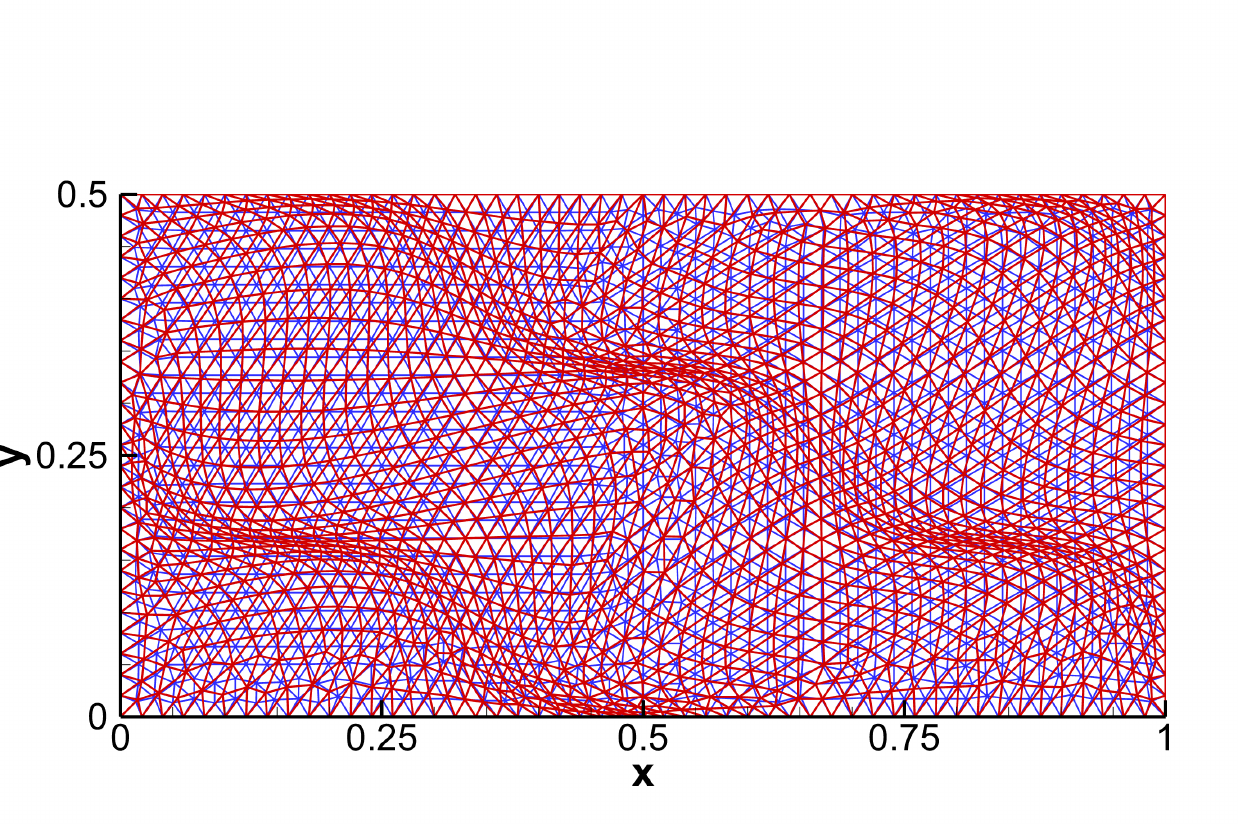}
\includegraphics[width=0.475\textwidth]{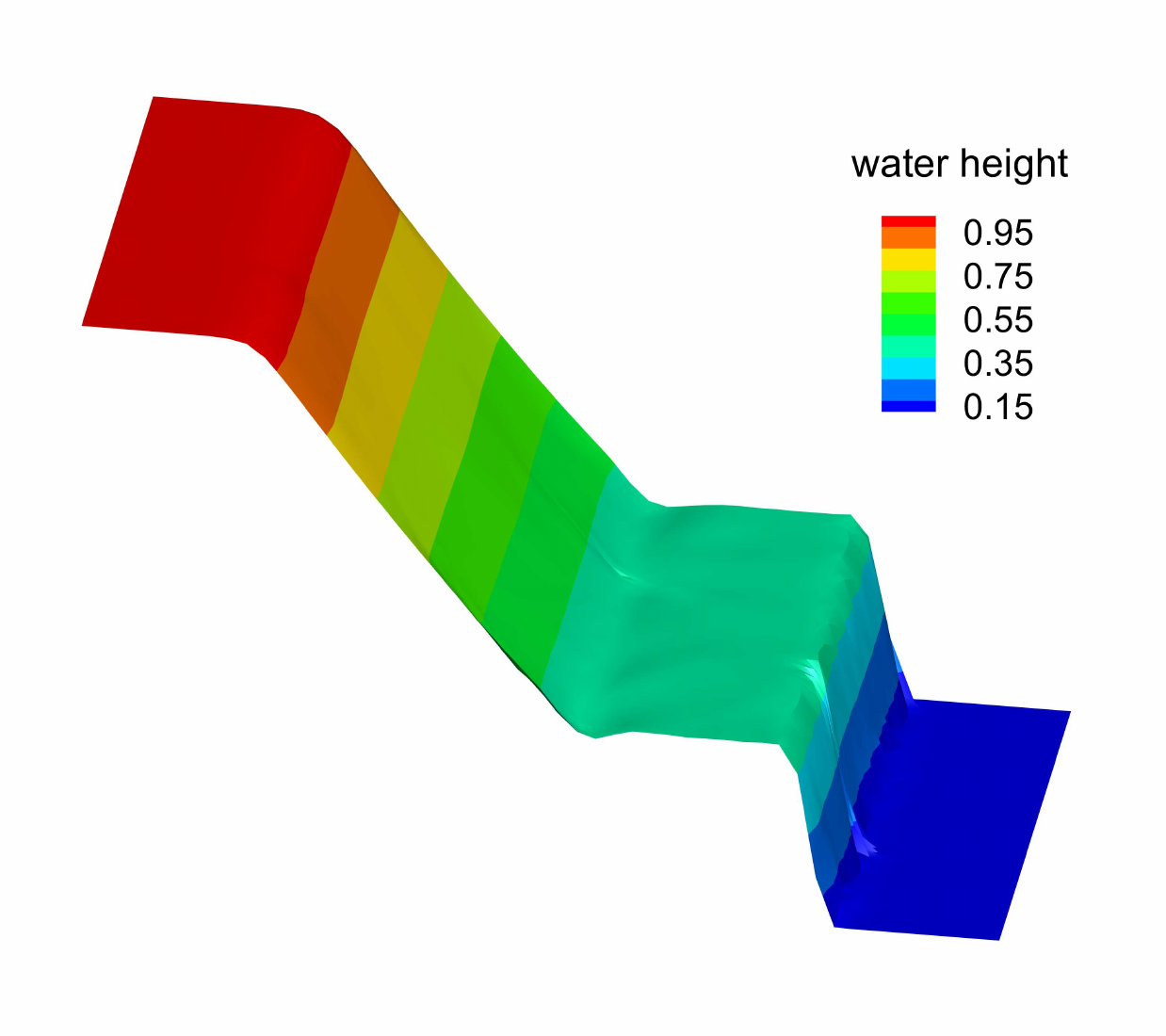}\\
\includegraphics[width=0.475\textwidth]{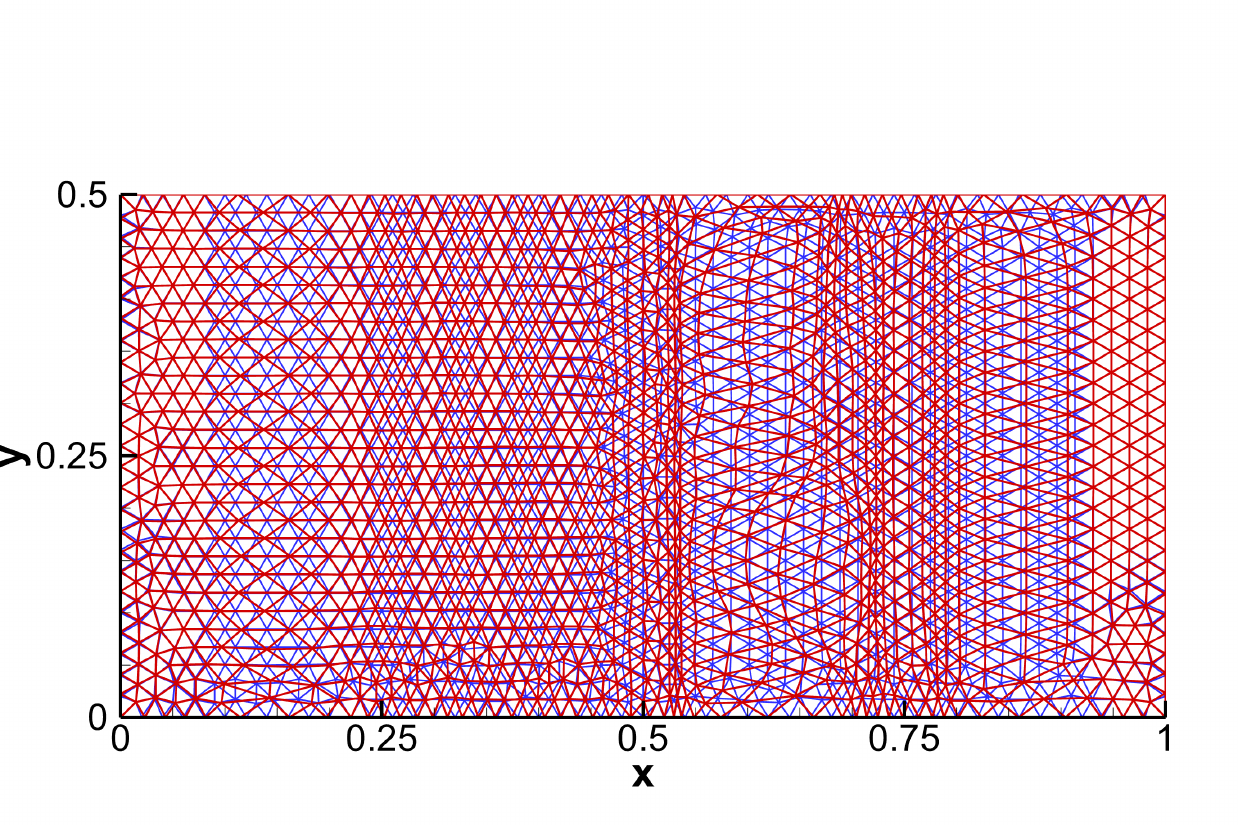}
\includegraphics[width=0.475\textwidth]{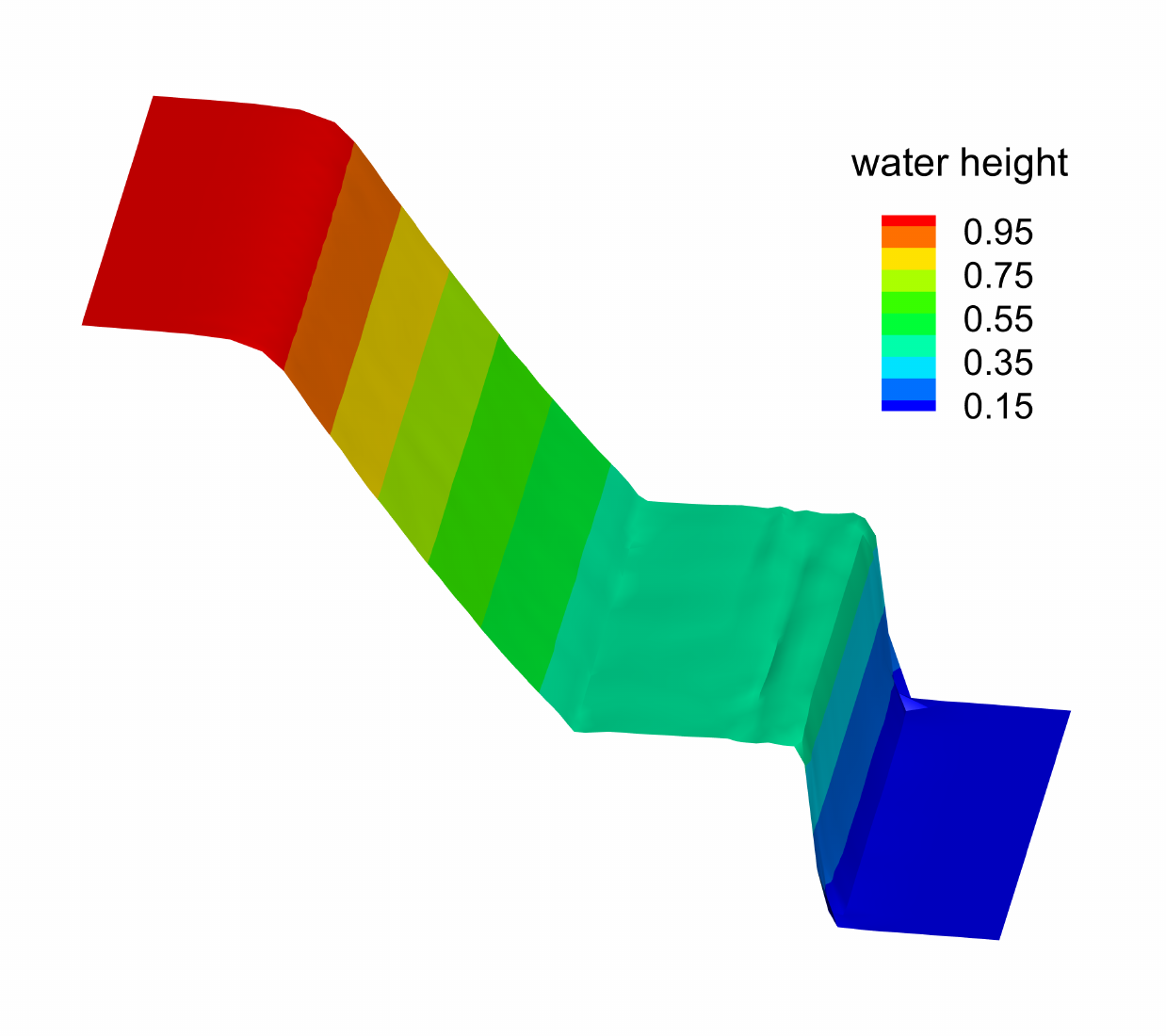}
\caption{\label{1d-dam-break-1} 1-D dam-break problem. Left: Comparison of the initial mesh (blue) with the final configuration (red) at $t=0.3$. Right: 3-D view of the final water height distribution.}
\end{figure}

\begin{figure}[!htb]
\centering
\includegraphics[width=0.495\textwidth]{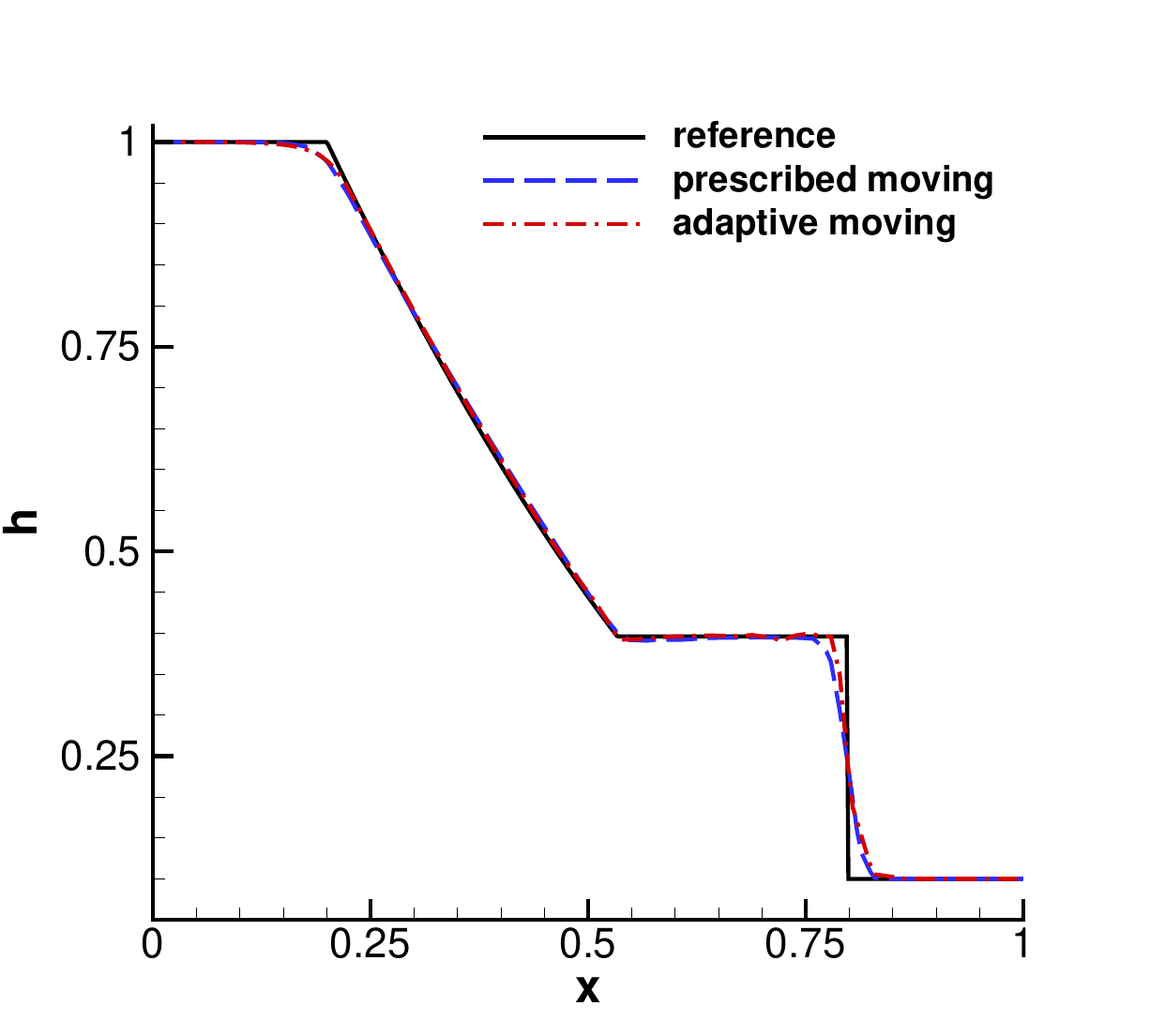}
\includegraphics[width=0.495\textwidth]{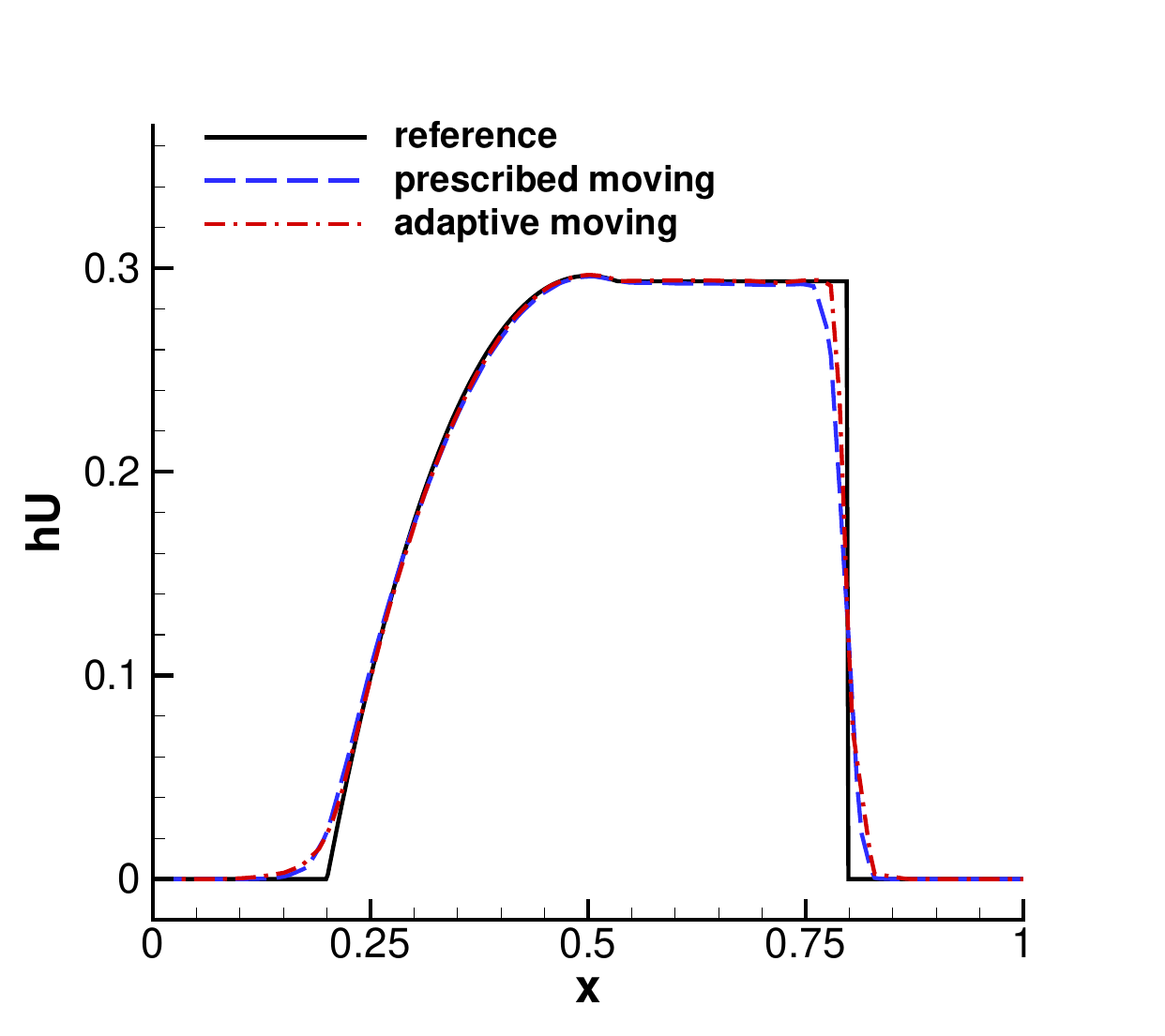}
\caption{\label{1d-dam-break-2} 1-D dam-break problem. Distributions of water height (left) and momentum (right) along the horizontal centerline of the domain by the compact GKS with different moving strategies.}
\end{figure}

\subsection{1-D dam-break problem}

To assess the performance of the proposed compact GKS on moving meshes in capturing strong discontinuities, we simulate the 1-D dam-break problem over a flat bottom. The initial condition is defined as
\begin{equation*}
(h,U,V) = \begin{cases}
(1,0,0),  ~~ & 0\leq x<0.5,\\
(0.1,0,0),   & 0.5\leq x\leq1.
\end{cases}
\end{equation*}
The simulation is conducted on a computational domain of $[0,1]\times[0,0.5]$, initialized with a uniform triangular mesh having a characteristic cell size of $\Delta X=0.02$.
The simulation is run until a final time of $t=0.3$.
Two distinct mesh motion strategies are employed. The first is a prescribed nodal motion, where the coordinates of the moving mesh are defined by
\begin{align*}
\mathbf{x}=\mathbf{x}_0 +0.05\mathrm{sin}(\pi t)\mathrm{sin}(3\pi x)\mathrm{sin}(6\pi y)\mathbf{e}.
\end{align*}
The second is an adaptive mesh motion method driven by variations in the flow variables, as presented in Section 5.

Fig. \ref{1d-dam-break-1} illustrates the results for both strategies. It compares the initial mesh (blue) with the final mesh (red) at $t=0.3$ and displays the corresponding water height distributions. For a more quantitative comparison, Figure \ref{1d-dam-break-2} plots the water height profiles along the horizontal centerline of the domain.
The results demonstrate that the proposed scheme yields accurate solutions for both prescribed and adaptive mesh motion. Notably, the solutions are free of numerical oscillations, effectively capturing the sharp discontinuity at the shock wave.

\subsection{2-D circular dam-break problem}
In this section, the dam break problem of a circular dam is simulated on a moving mesh.
The initial condition of the test case is given by
\begin{equation*}
h = \begin{cases}
1.0, ~~ (x-5)^2+(y-5)^2 <4.0,\\
0.5, ~~ \mathrm{otherwise}.
\end{cases}
\end{equation*}
The initial velocity is $(U,V)=(0,0)$ in the computational domain $[0,10]\times[0,10]$.
The wall boundary condition is adopted on all boundaries. An initially uniform triangular mesh with a cell size of $\Delta X=0.15$ is used in the computation.

\begin{figure}[!htb]
\centering
\includegraphics[width=0.32\textwidth]{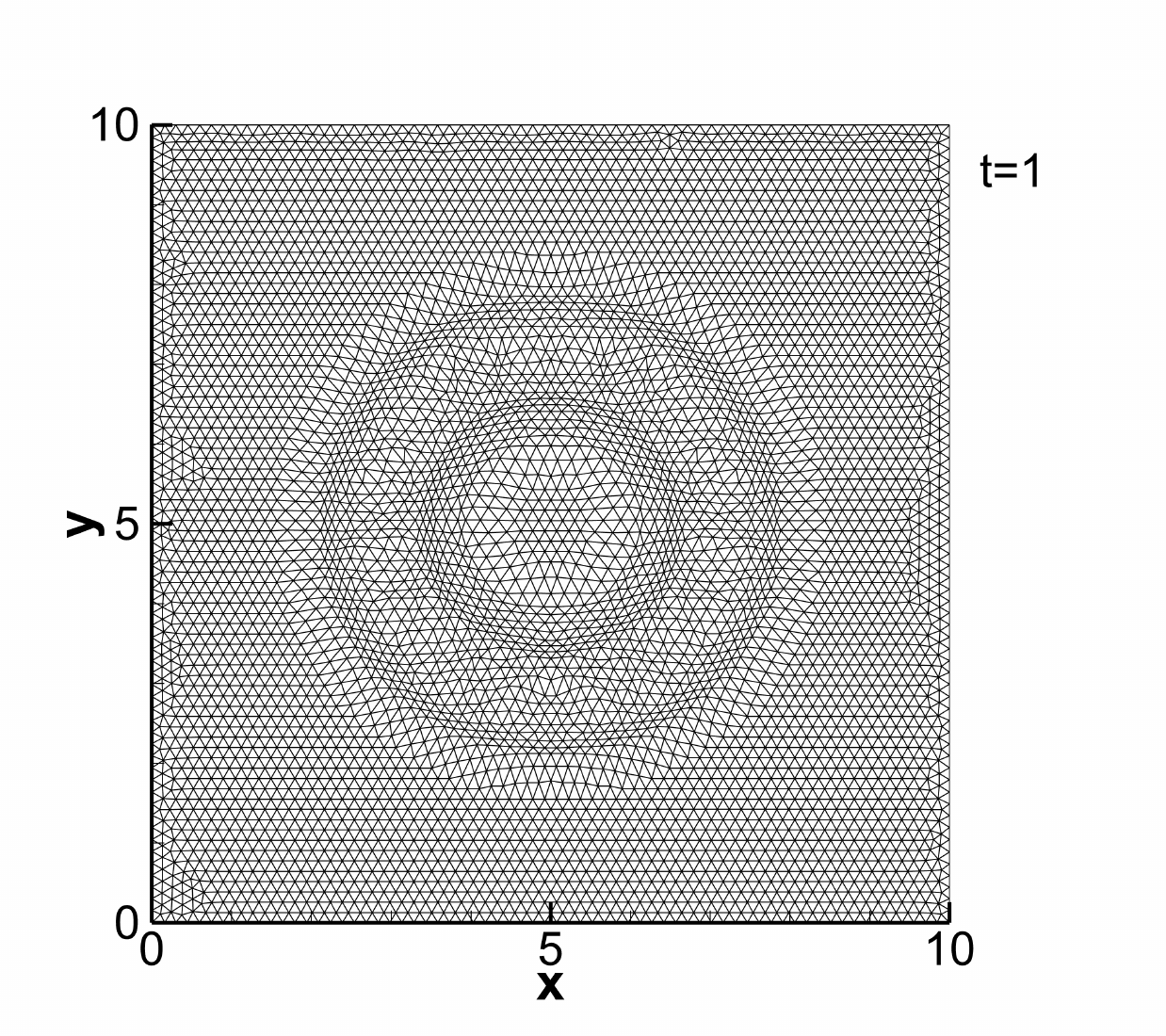}
\includegraphics[width=0.32\textwidth]{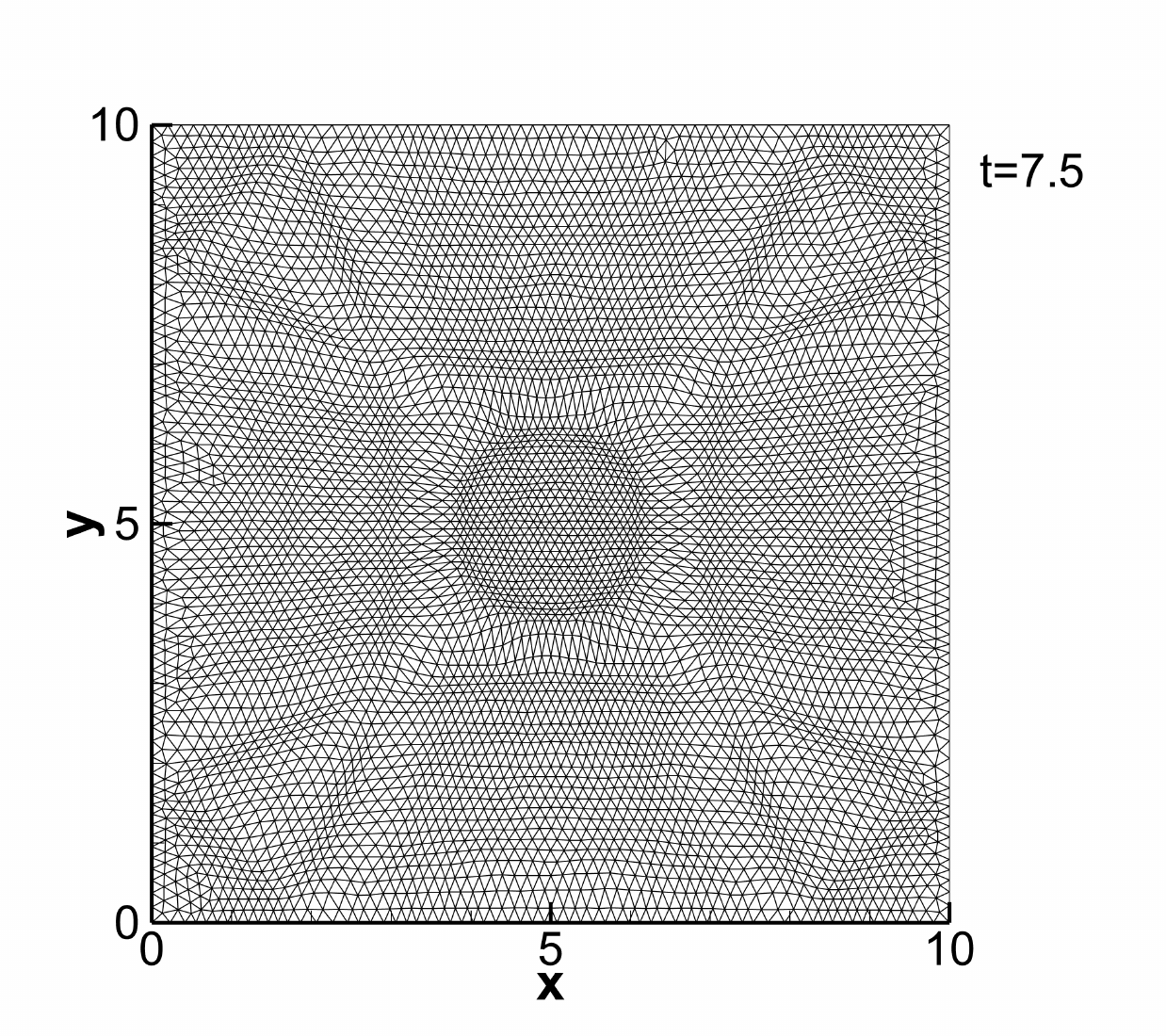}
\includegraphics[width=0.32\textwidth]{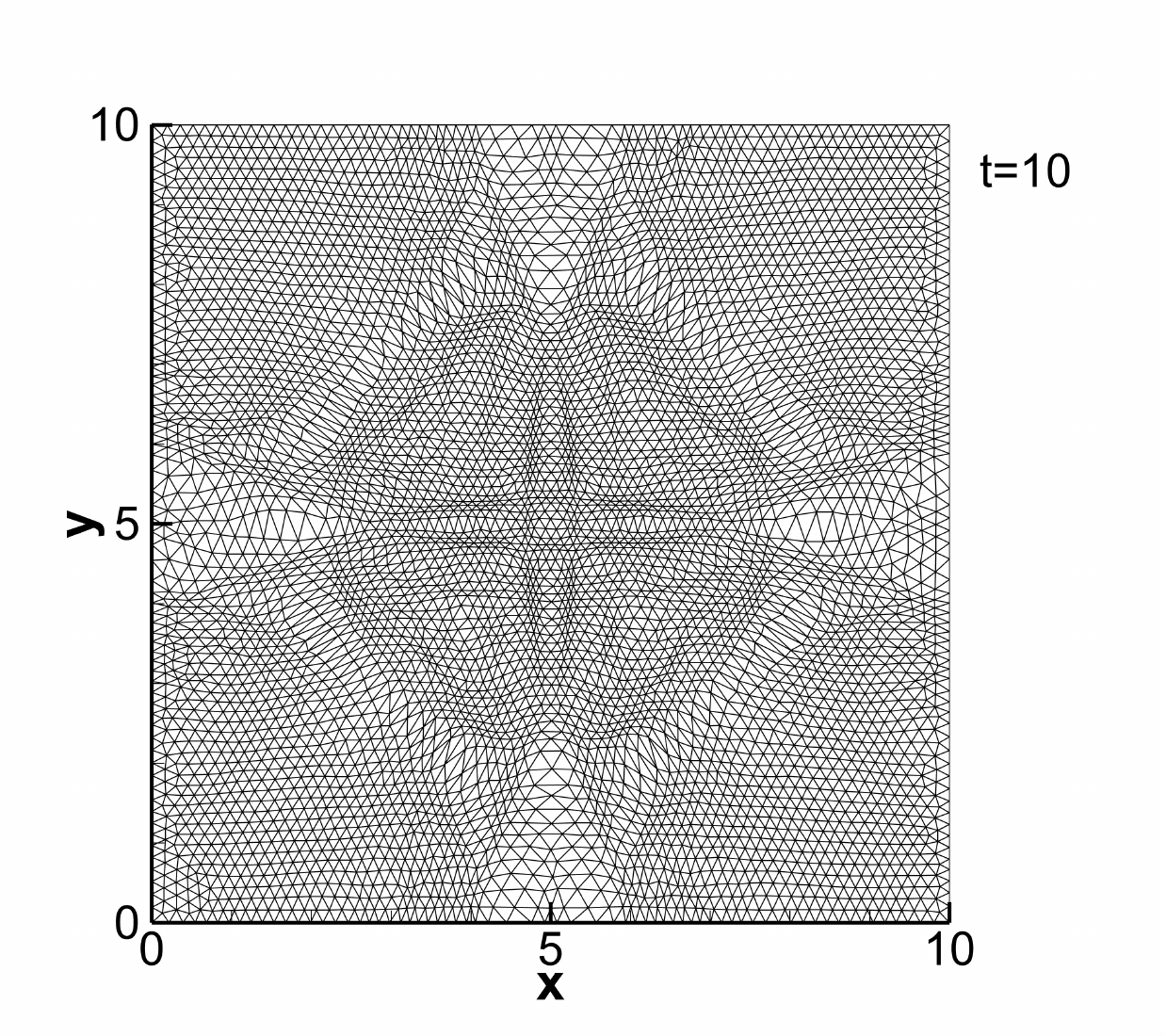} \\
\includegraphics[width=0.32\textwidth]{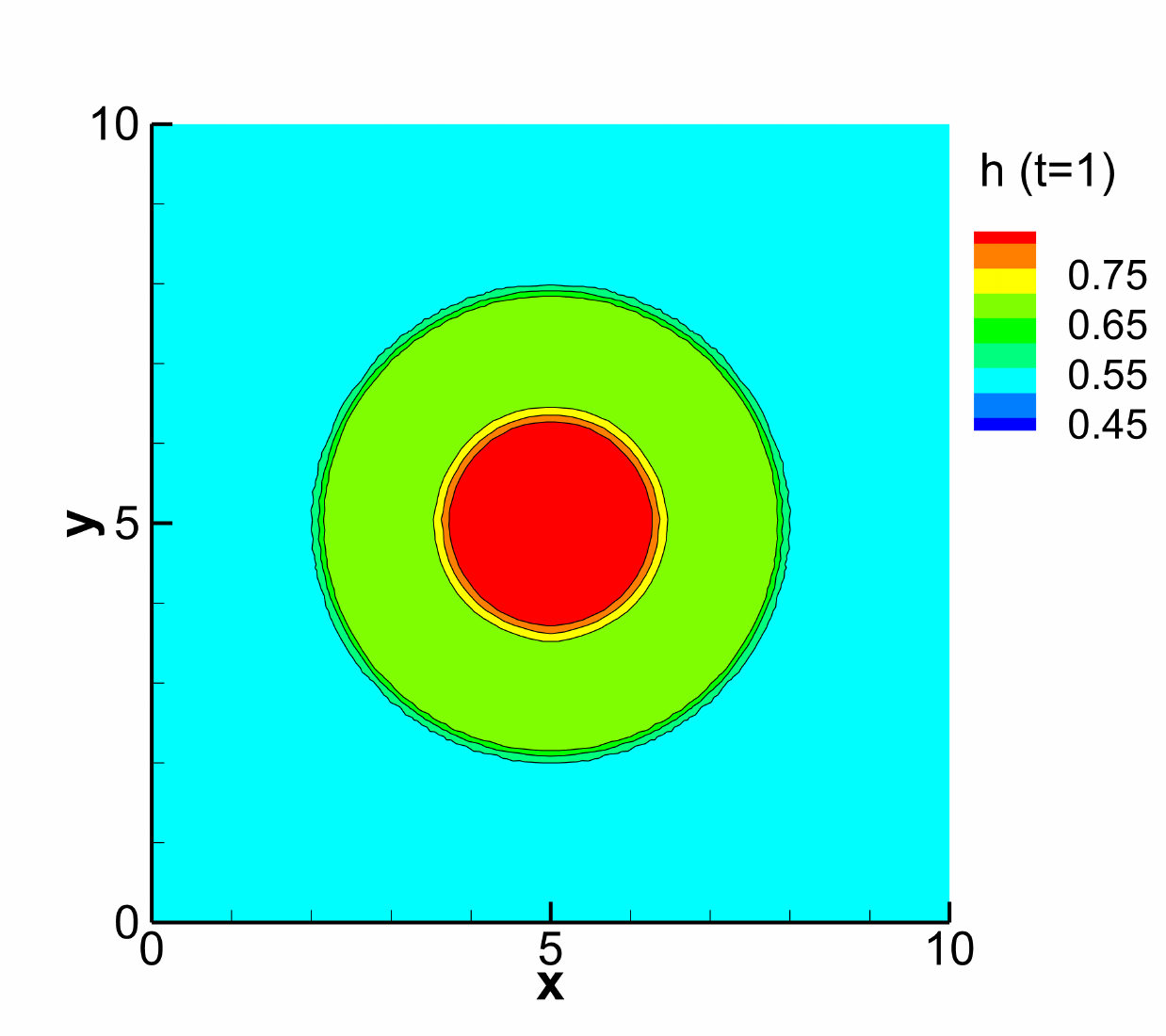}
\includegraphics[width=0.32\textwidth]{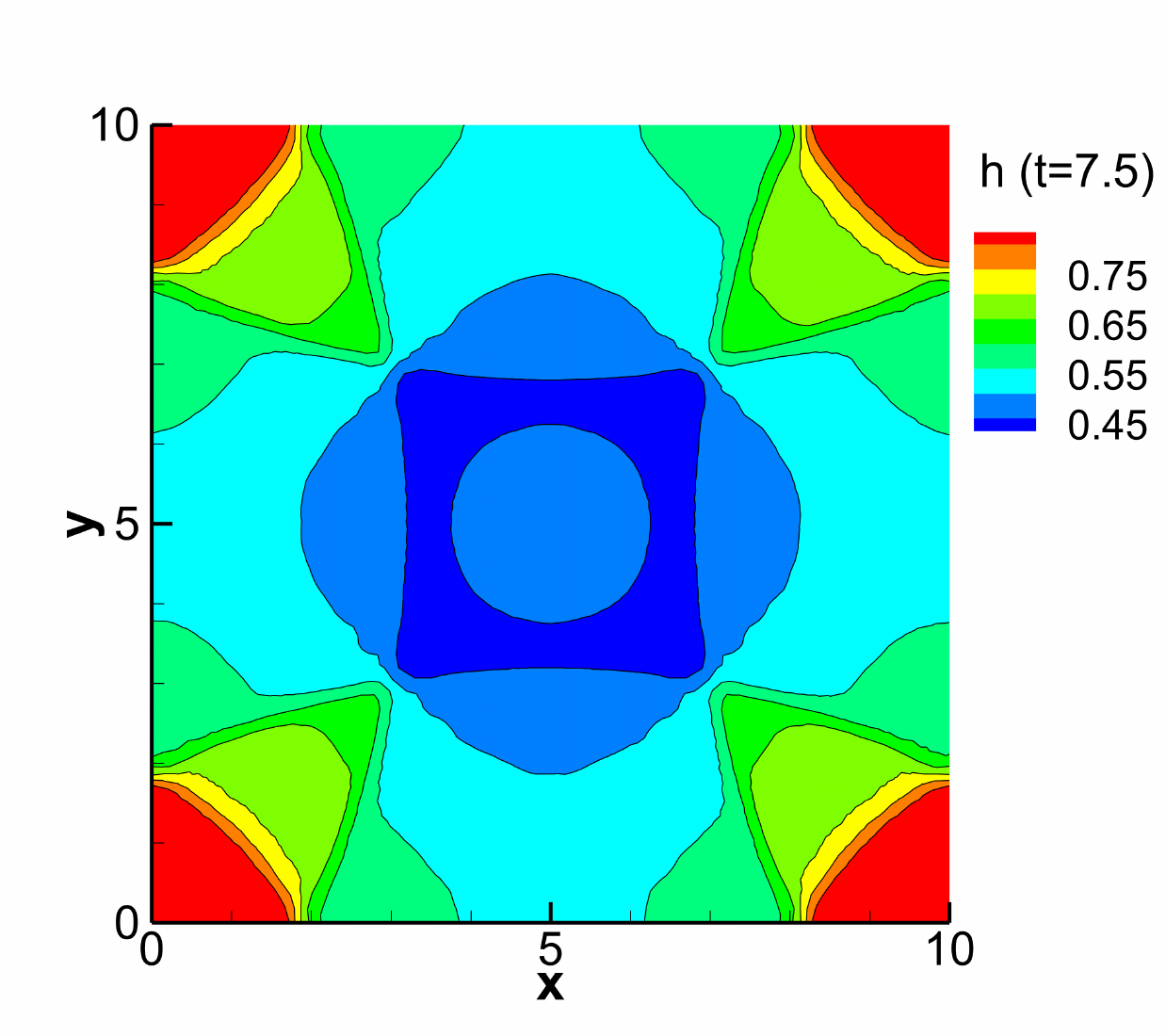}
\includegraphics[width=0.32\textwidth]{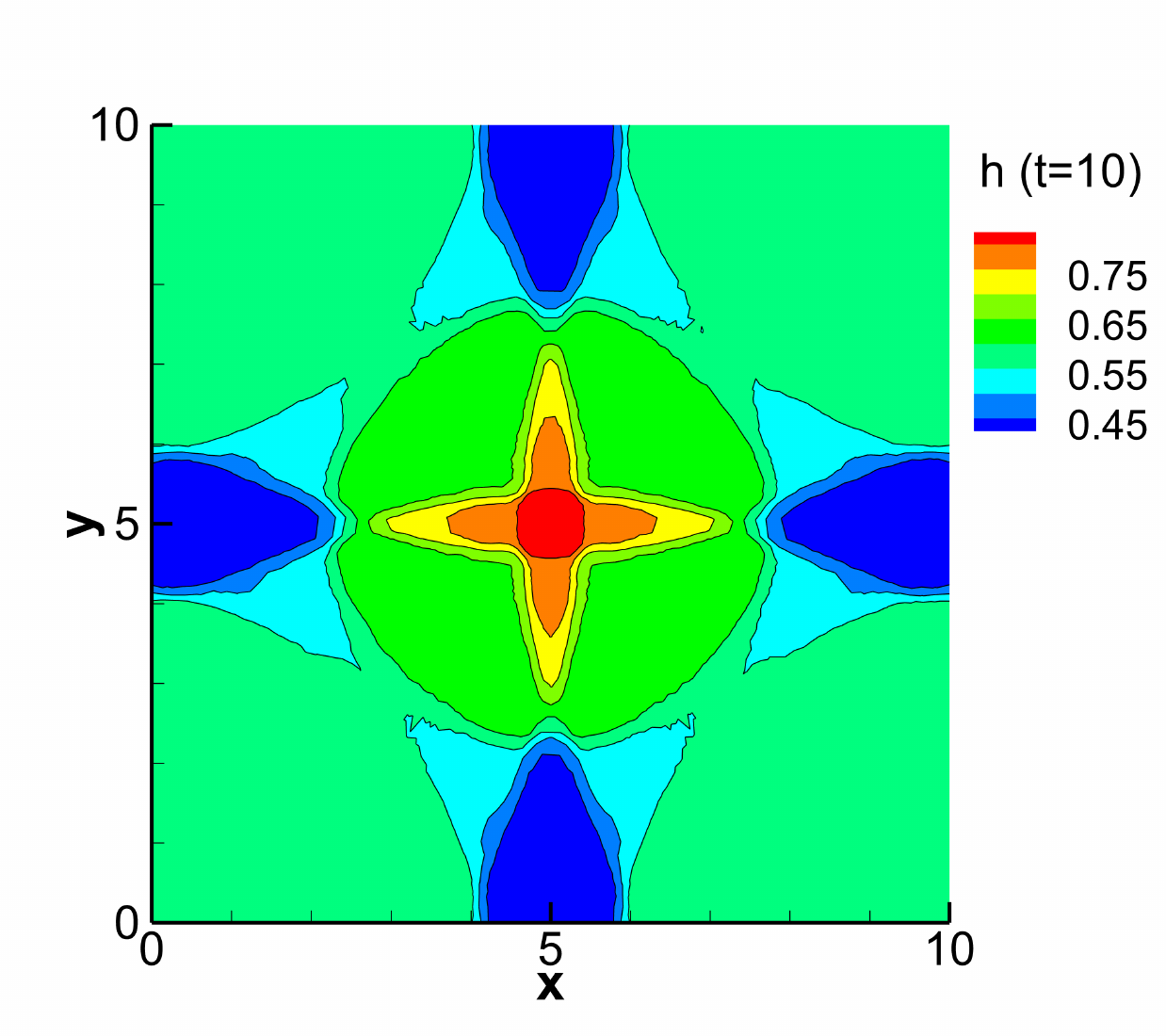}
\caption{\label{2d-circular-dam-break} 2-D circular dam-break problem: the adaptively moved mesh (up) and the water height contours (down) at different simulation times $t=1.0$, $t=7.5$ and $t=10.0$.}
\end{figure}

The mesh after adaptive moving and the distribution of water surface height at $t=1.0$, $t=7.5$ and $t=10.0$ are presented in Fig. \ref{2d-circular-dam-break}.
An initially cylindrical water column generates axisymmetric shock waves that propagate outward and reflect at the square wall. The reflected shock waves converge towards the center, leading to complex symmetric structures on the water surface.
The present ALE compact GKS captures the complex wave patterns with high resolution, preserves the symmetric structure of the numerical solution, and exhibits no spurious oscillations. Comparison of the mesh and water height distributions indicates that the mesh undergoes adaptive refinement during the computation.

\subsection{2-D dam-break in an irregular domain}
The 2-D dam-break problem with an irregular domain in \cite{dambreak-1998,dambreak-2007} is used in the current study to validate the compact GKS on an adaptive moving unstructured mesh.
The computational domain is set as follows. The length of the dam breach is $75$ and it starts at $y=95$.
The dam itself has a width of $10$ and its left side is located at $x = 95$.
The bottom topography is flat.
At $t=0$ the stationary water surface has a discontinuity with $h_l=10$ and $h_r=5$ across the breach.
The water level is set as
\begin{equation*}
h = \begin{cases}
10,  ~0\leq x<95,\\
5,   ~~95\leq x.
\end{cases}
\end{equation*}
The boundary condition on the far right is the free boundary, and the other boundary conditions are the non-penetration slip wall boundaries.
The gravitational acceleration is $G=9.812$.

\begin{figure}[!htb]
\centering
\includegraphics[width=0.495\textwidth]{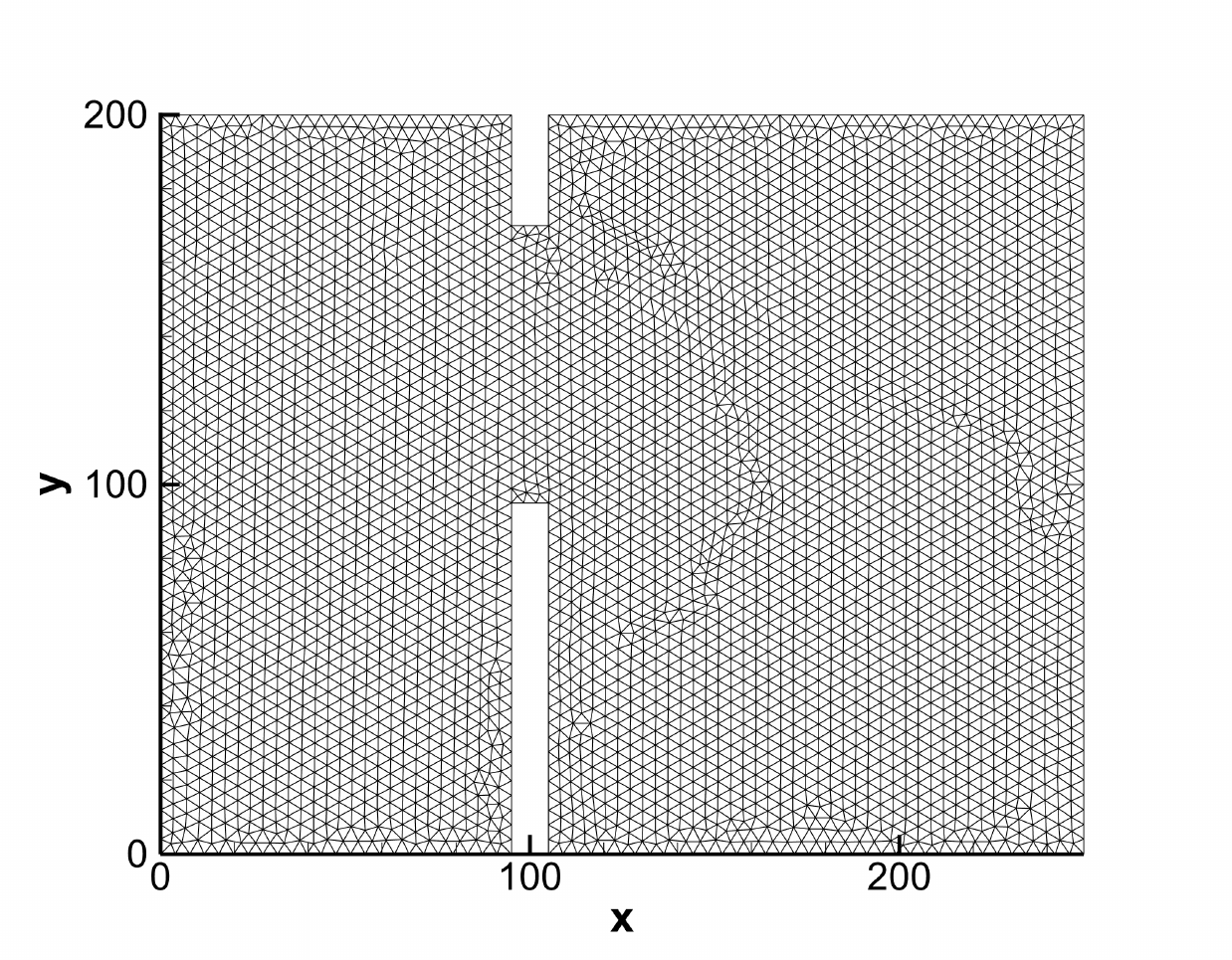}
\includegraphics[width=0.495\textwidth]{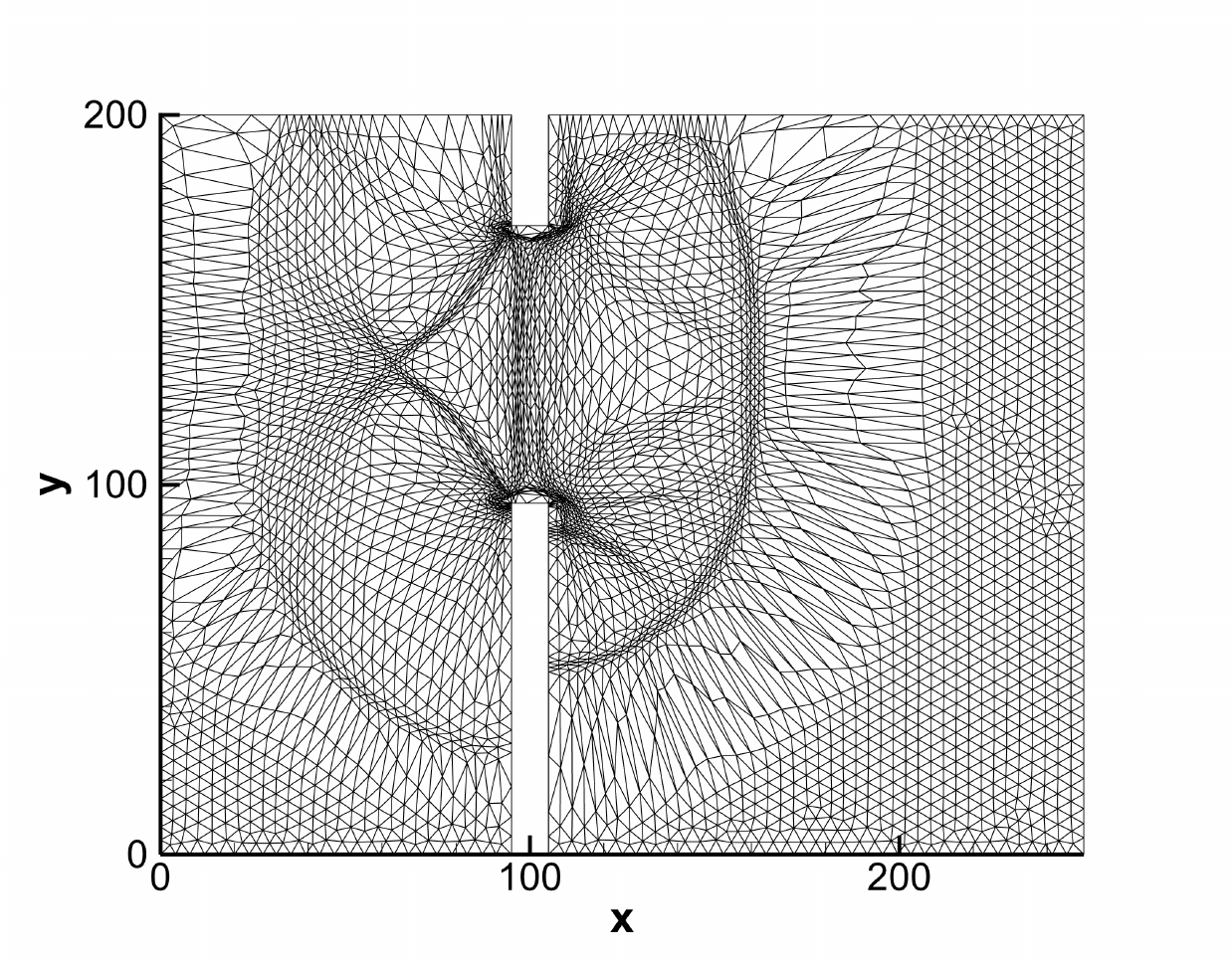}
\caption{\label{2d-dam-break-1} 2-D dam-break in an irregular domain: the initial mesh (left) and the adaptively moved mesh (right) at the end of the simulation at time $t=7.2$.}
\end{figure}

\begin{figure}[!htb]
\centering
\includegraphics[width=0.495\textwidth]{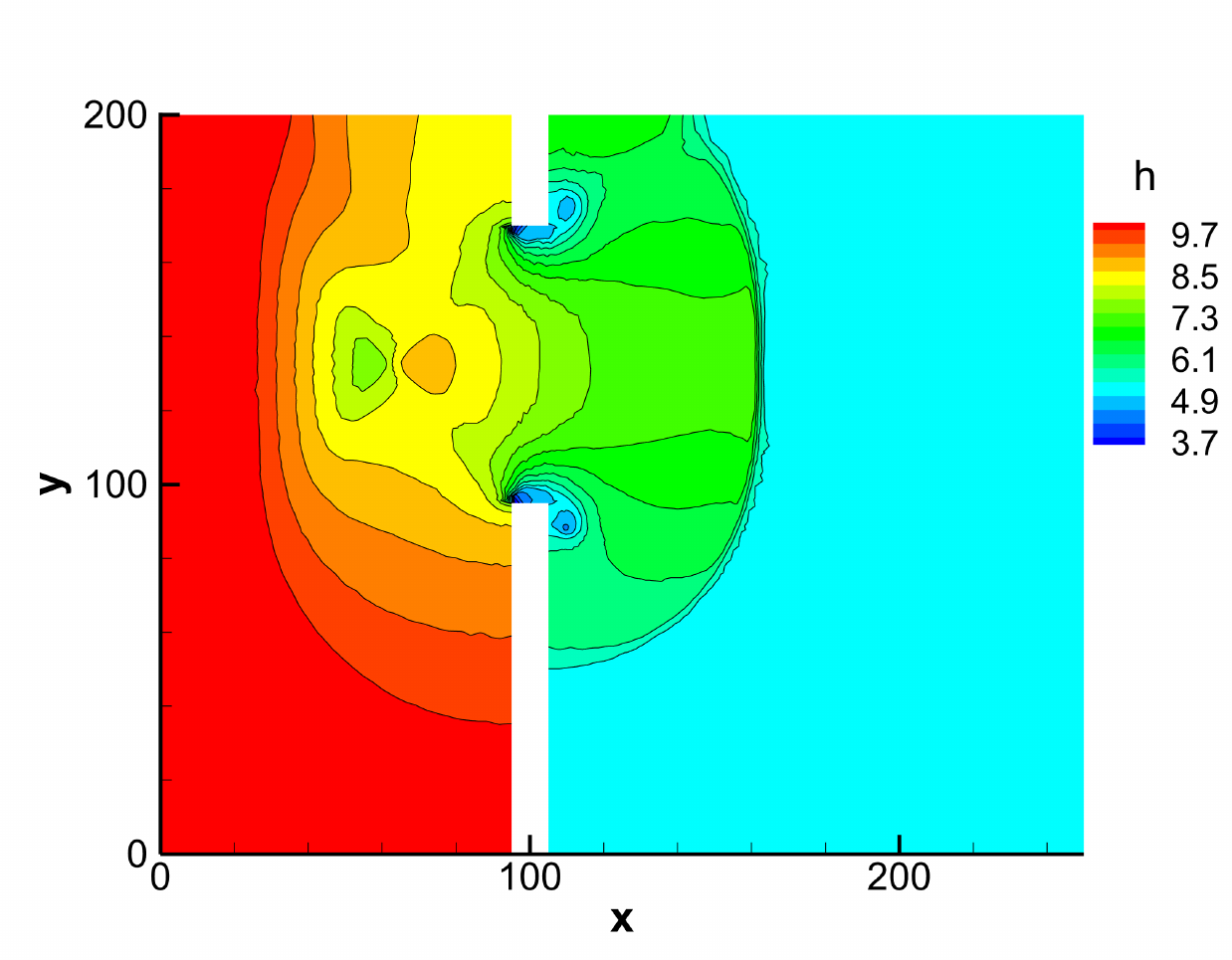}
\includegraphics[width=0.495\textwidth]{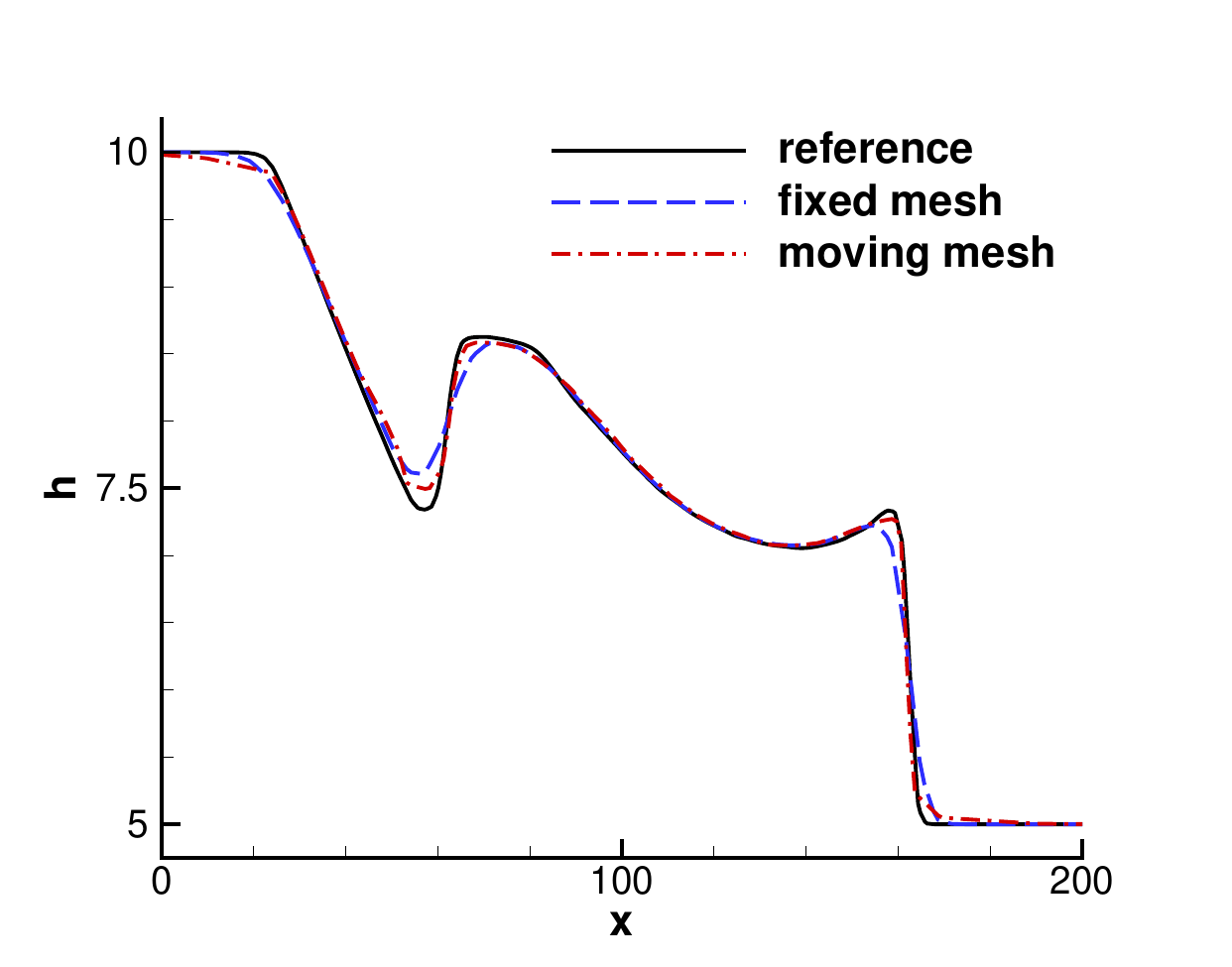}
\caption{\label{2d-dam-break-2} 2-D dam-break in an irregular domain: the 2-D water height distribution (left) and the distribution of water height along the horizontal centerline of the dam (right).}
\end{figure}

The initial mesh spacing is $\Delta X=4.0$, and the initial mesh is shown in Fig. \ref{2d-dam-break-1}. At $t=7.2$, the adaptively moved mesh is displayed in the right panel of Fig. \ref{2d-dam-break-1}. Comparing this with the water height distribution in the left panel of Fig. \ref{2d-dam-break-2} shows that the mesh concentrates in regions with larger variations in water height.
The right panel of Fig. \ref{2d-dam-break-2} plots the water height along the horizontal centerline of the dam. The reference solution is computed on a refined fixed mesh using the fourth-order compact GKS \cite{zhao2021-swe}. The two solutions agree closely, while the moving-mesh ALE simulation uses approximately one quarter of the degrees of freedom of the refined-mesh computation. In the refined-mesh run, the local cell size near the dam is reduced to $1/13$ of the initial adaptive mesh spacing.

\section{Conclusion}

In this study, we develop a space-time coupled ALE compact GKS for solving the SWE on moving triangular meshes. The temporal update is performed directly on the physical moving mesh, obviating data remapping.
The proposed scheme rigorously satisfies the GCL for arbitrarily moving meshes and preserves the well-balanced property over non-flat bottom topography. These properties are verified through both theoretical analysis and numerical experiments. Moreover, the scheme employs a fourth-order compact spatial reconstruction, achieving non-oscillatory and robust solutions even in the presence of discontinuities.

The high accuracy and strong robustness of the proposed scheme arise from several key design features. First, the numerical fluxes are directly derived by accounting for the effects of mesh motion, including contributions from the spatial temporal nonuniformity of the flow field and the convective fluxes associated with the swept areas of moving cell interfaces. Second, to consistently represent the evolving bottom topography on moving meshes, a governing equation for topography is formulated and discretized using a compatible space-time coupled approach. Third, the compact GKS provides time-accurate evolution of both flow variables and fluxes for the ALE flux computation, which simultaneously serves as the foundation for the high-order compact reconstruction.

The proposed space-time coupled ALE scheme is particularly suitable for problems involving adaptive mesh refinement to resolve local small-scale flow structures and fluid structure interaction phenomena. Its broader practical applications will be further explored in future work. In summary, this study presents a high-order ALE compact GKS for the SWE with non-flat bottom topography within a unified space-time discretization framework. The construction methodology can be naturally extended to general hyperbolic conservation laws without source terms.

\section*{Acknowledgments}
This research was supported by the Research Grants Council Areas of Excellence (AoE) Scheme (AoE/P-601/23N-D - MATH), and by CORE as a joint research center for ocean
research between Laoshan Laboratory and HKUST.

\section*{Appendix}

The space-time ALE compact GKS in this study requires, at each cell interface, the flow variables and fluxes, their temporal derivatives, and the spatial gradients of the flow variables. All of these quantities are obtained from the time-dependent solution of the gas distribution function in the GKS, and the corresponding formulas are presented here. Further details of the GKS for the SWE can be found in \cite{xu2002-swe,zhao2021-swe}.

The second-order accurate evolution solution of the gas distribution function is given as
\begin{equation}\label{SWE-2nd-order}
\begin{split}
f(t,\textbf{u})&=\overline{g}(\mathbf{u})\big[ C_1+ C_2 \big( \overline{\mathbf{a}}^l \cdot\mathbf{u}H(u) +\overline{\mathbf{a}}^r \cdot\mathbf{u}(1-H(u)) \big) +C_3\overline{A} \big]\\
                            &+C_2\overline{g}(\mathbf{u}) \big[-2 \alpha_{k,m} \overline{\lambda} \big( \nabla\Phi^l H(u)+\nabla\Phi^r(1-H(u)) \big) \cdot (\mathbf{u}-\overline{\mathbf{U}}) \big]\\
                            &+C_4\big[g^l(\mathbf{u})H(u)+g^r(\mathbf{x},0,\mathbf{u})(1-H(u))\big] \\
                            &+C_5g^l(\mathbf{u})\big[\mathbf{a}^l\cdot \mathbf{u} -2 \alpha_{k,m} \lambda^l \nabla\Phi^l \cdot(\mathbf{u}-\mathbf{U}^l) \big]H(u) \\
                            &+C_5g^r(\mathbf{u})\big[\mathbf{a}^r\cdot \mathbf{u} -2 \alpha_{k,m} \lambda^r \nabla\Phi^r \cdot(\mathbf{u}-\mathbf{U}^r) \big](1-H(u)),
\end{split}
\end{equation}
where $f$ is the gas distribution function, $\textbf{u}=(u,v)$ is the particle velocity, $\nabla\Phi^{l,r}$ are determined by the left and right values of the source term at the interface, and $C_k$ are the coefficients related to $t$.
$\overline{g}$ and $g^{l,r}$ are the equilibrium states at $t^n$ determined by $W^e(t^n)$ and $W^{l,r}(t^n)$, respectively. $\overline{\mathbf{a}}^{l,r}$ and $\mathbf{a}^{l,r}$ are the space derivatives of $\overline{g}$ and $g^{l,r}$, and $\overline{A}$ represents the corresponding time derivative of $\overline{g}$.
The equilibrium state $g$ is a Maxwellian distribution function,
\begin{equation*}
\begin{split}
g=h\big(\frac{\lambda}{\pi}\big)e^{-\lambda(\mathbf{u}-\mathbf{U})^2},
\end{split}
\end{equation*}
where $\lambda$ is defined by $\lambda=1/Gh$.
$\alpha_{k,m}~(k=1,2,~m=1,2,3)$ are constants introduced to ensure the well-balanced evolution solution, with $(\alpha_{1,1},\alpha_{1,2},\alpha_{1,3})=(1,3/4,1/4)$ and $\alpha_{2,m}=1$. The index $k$ indicates that $\alpha_{k,m}$ corresponds to the $k$-th component of the vector $\nabla\Phi^{l,r}$, and $m$ corresponds to the order of the moment with respect to the microscopic velocity $u$, i.e., $<u^m>$.
The moment operation $<\cdots>$ is defined as
\begin{equation*}
<\cdots>=\int (\cdots)f \mathrm{d}u\mathrm{d}v.
\end{equation*}
For the second-order evolution solution, the above $f(t)$ can be approximated through a linearization in time \cite{zhaocompact_tri}
\begin{equation*}
\hat{f}(t)=f^n+t f_t^n.
\end{equation*}
The two coefficients $f^n$ and $f_t^n$ are calculated as follows
\begin{align*}
    f^n&=\big(4\bar{f}(\Delta t/2) - \bar{f}(\Delta t)\big)/\Delta t,\\
    f_t^n&=4\big(\bar{f}(\Delta t) - 2\bar{f}(\Delta t/2)\big)/{\Delta t}^2,
\end{align*}
where $\bar{f}(\Delta t)$ and $\bar{f}(\Delta t/2)$ are the time integrations of $f(t)$ over the interval $[t_n, t_n + \Delta t]$ and $[t_n, t_n + \Delta t/2]$, respectively.

The numerical fluxes and their time derivatives can be obtained by taking moments of $\hat{f}(t)$ and $\hat{f}_t(t)$ at $t = t_n$
\begin{equation*}
\boldsymbol{F}^n = \int u f^n \boldsymbol{\psi} \, \text{d}u\text{d}v,\quad
\boldsymbol{F}_t^n = \int u f_t^n \boldsymbol{\psi} \, \text{d}u\text{d}v.
\end{equation*}
Simultaneously, the flow variables and their time derivatives can be given by
\begin{equation*}
\boldsymbol{W}^n = \int f^n \boldsymbol{\psi} \, \text{d}u\text{d}v,\quad
\boldsymbol{W}_t^n = \int f_t^n \boldsymbol{\psi} \, \text{d}u\text{d}v.
\end{equation*}
While interface $\mathbf{W}(t)$ can be discontinuous, and the modeled discontinuous interface evolution solutions are given as
\begin{equation*}
\begin{split}
{\bf W}^l(t^{n+1})&= (1-e^{-\Delta t/\tau_0}){\overline{\bf W}}(t^{n+1}) +e^{-\Delta t/\tau_0} {\bf W}_0^l(t^{n+1}), \\
{\bf W}^r(t^{n+1})&= (1-e^{-\Delta t/\tau_0}){\overline{\bf W}}(t^{n+1}) +e^{-\Delta t/\tau_0} {\bf W}_0^r(t^{n+1}),
\end{split}
\end{equation*}
where $\tau_0=C \big|\frac{h_l^2-h_r^2}{h_l^2+h_r^2}\big|\Delta t$ with the constant $C=5$. $\overline{\bf W}(t^{n+1})$ and ${\bf W}_0^{l,r}(t^{n+1})$ are determined by the corresponding $\boldsymbol{W}^n$ and $\boldsymbol{W}_t^n$.

In addition, the initial left and right interface states, $\mathbf{W}^{l,r}$, together with their spatial derivatives, $\partial \mathbf{W}^{l,r}/\partial s$ ($s \in {x,y}$), are obtained via spatial reconstruction and used to determine $g^{l,r}$. The initial interface equilibrium state, $\overline{\mathbf{W}}$, is then specified by the conservation constraints as
\begin{equation*}
\overline{\mathbf{W}}= \int \overline{g} \boldsymbol{\psi}\text{d}u\text{d}v = \int_{u>0} g^l \boldsymbol{\psi}\text{d}u\text{d}v +\int_{u<0} g^r \boldsymbol{\psi} \text{d}u\text{d}v,
\end{equation*}
and the initial interface equilibrium spatial derivatives $\partial \overline{\mathbf{W}}/\partial s$ is modeled as
\begin{equation*}
\partial \overline{\mathbf{W}}/\partial s= \int \partial \overline{g}/\partial s \boldsymbol{\psi} \text{d}u\text{d}v = \int_{u>0} \partial g^l/\partial s \boldsymbol{\psi} \text{d}u\text{d}v +\int_{u<0} \partial g^r/\partial s \boldsymbol{\psi} \text{d}u\text{d}v.
\end{equation*}

\bibliographystyle{ieeetr}
\bibliography{compact-gks}

\end{document}